\documentclass[12pt]{amsart}
\usepackage{amsmath, amssymb, booktabs,mathtools}
\usepackage[margin=1in]{geometry}
\usepackage{color,xcolor}
\usepackage[allcolors=black,pdftex]{hyperref}
\usepackage[T1]{fontenc}
\usepackage{tikz}
\usetikzlibrary{calc}
\usepackage{comment}
\usepackage{longtable,array}

\newcolumntype{C}{>{$}c<{$}} 
\newcolumntype{L}{>{$}l<{$}} 
\newcolumntype{R}{>{$}r<{$}} 

\theoremstyle{plain}

\newtheorem{theorem}{Theorem}[section]
\newtheorem{lemma}[theorem]{Lemma}
\newtheorem{proposition}[theorem]{Proposition}
\newtheorem{corollary}[theorem]{Corollary}

\theoremstyle{definition}
\newtheorem{definition}[theorem]{Definition}
\newtheorem{example}[theorem]{Example}

\theoremstyle{remark}
\newtheorem{remark}{Remark}[section] 

\DeclareMathOperator{\n}{N}     

\DeclareMathOperator{\SL}{SL} 
\DeclareMathOperator{\Gal}{Gal} 

\newcommand{\GL}[1]{\operatorname{GL}_2(#1)} 
\newcommand{\CL}{\operatorname{Cl}_{K}} 
\newcommand{\CLno}{h_K} 
\newcommand{\disc}{\operatorname{disc}} 
\newcommand{\SO}{\operatorname{SO}} 
\newcommand{\SU}{\operatorname{SU}} 
\newcommand{\hy}{\mathbb{H}} 
\newcommand{\CC}{\mathbb{C}} 
\newcommand{\RR}{\mathbb{R}} 
\newcommand{\QQ}{\mathbb{Q}} 
\newcommand{\ZZ}{\mathbb{Z}} 
\newcommand{\PP}{\mathbb{P}} 
\renewcommand{\AA}{\mathbb{A}} 
\newcommand{\TT}{\mathbb{T}} 
\newcommand{\OO}{\mathcal{O}_K} 

\newcommand{\go}{\GL{\OO}}
\newcommand{\ur}{\text{ur}}
\newcommand{\new}{\text{new}}

\newcommand{\mf}[1]{\mathfrak{#1}}
\newcommand{\mat}[1]{\begin{bmatrix} #1 \end{bmatrix}}
\newcommand{\leg}[2]{\left(\frac{#1}{#2}\right)}   

\newcommand{\pp}{\mathfrak{p}} 
\renewcommand{\aa}{\mathfrak{a}} 
\newcommand{\bb}{\mathfrak{b}}
\newcommand{\dd}{\mathfrak{d}} 
\newcommand{\qq}{\mathfrak{q}} 
\newcommand{\nn}{\mathfrak{n}}
\newcommand{\mm}{\mathfrak{m}}
 
\DeclarePairedDelimiterX\set[1]{\lbrace}{\rbrace}{\def\given{\;\delimsize\vert\;}#1}
\newcommand{\Qroot}[1]{\QQ(\sqrt{#1})}

\newcommand{\software}[1]{\textsc{#1}{}}

\newcommand{\Magma}{\software{Magma}}
\newcommand{\Github}{\software{GitHub}}

\newcommand{\lmfdb}{\href{https://www.lmfdb.org}{LMFDB}}
\newcommand{\lmfdbfield}[1]{\href{https://www.lmfdb.org/NumberField/#1}{\textsf{#1}}}
\newcommand{\lmfdbbmftop}{\href{https://www.lmfdb.org/ModularForm/GL2/ImaginaryQuadratic}{LMFDB}}

\newcommand{\lmfdbbmf}[3]{\href{https://www.lmfdb.org/ModularForm/GL2/ImaginaryQuadratic/#1/#2/#3}{\textsf{#1-#2-#3}}}

\newcommand{\lmfdbcmf}[4]{\href{https://beta.lmfdb.org/ModularForm/GL2/Q/holomorphic/#1/#2/#3/#4/}{\textsf{#1.#2.#3.#4}}}
\newcommand{\lmfdbec}[3]{\href{https://www.lmfdb.org/EllipticCurve/Q/#1/#2/#3}{\textsf{#1.#2#3}}}

\newcommand{\lmfdbecnf}[4]{\href{https://www.lmfdb.org/EllipticCurve/#1/#2/#3}{\textsf{#1.#2-#3.#4}}}

\begin{document}
\title[Bianchi modular forms]{Bianchi Modular Forms over Imaginary Quadratic Fields with arbitrary class group}

\thanks{Last updated: \today}
\author[Cremona]{John Cremona}
\address{Warwick Mathematics Institute\\
  University of Warwick\\
  Coventry, CV4 7AL, UK}
\email{J.E.Cremona@warwick.ac.uk}
\urladdr{\url{https://johncremona.github.io/}}
\author[Thalagoda]{Kalani Thalagoda}
\address{Department of Mathematics and Statistics\\ 
Tulane University\\New Orleans, LA 70118}
\email{kthalagoda@tulane.edu}
\urladdr{\url{https://kalani-thalagoda.com}}
\author[Yasaki]{Dan Yasaki}
\address{Department of Mathematics and Statistics\\ 
University of North Carolina at Greensboro\\Greensboro, NC 27412}
\email{d\_yasaki@uncg.edu}
\urladdr{\url{https://go.uncg.edu/d_yasaki}}


\begin{abstract}
Let $K$ be an imaginary quadratic field and let $\mathcal{O}_K$ be its ring of integers. For an integral ideal $\nn$ of $\OO$, let $\Gamma_0({\nn})$ be the congruence subgroup of level ${\nn}$ consisting of matrices in $\GL{\OO}$ that are upper triangular mod ${\nn}$. In this paper, we discuss techniques to compute the space of Bianchi modular forms of level $\Gamma_0({\nn})$ as a Hecke module in the case where $K$ has arbitrary class group.  Our algorithms and computations extend and complement that carried out for fields of class number $1$, $2$, and $3$ by the first author, and by his students Bygott and Lingham in unpublished theses.  We give details and several examples for $K=\QQ(\sqrt{-17})$, whose class group is cyclic of order $4$, including a proof of modularity of an elliptic curve over this field. We also give an overview of the results obtained for a wide range of imaginary quadratic fields, which are tabulated in the L-functions and modular forms database (\href{https://www.lmfdb.org/}{LMFDB}).

\end{abstract}

\maketitle
\section{Introduction}
The goal of this paper is to study Bianchi modular forms over arbitrary imaginary quadratic fields from a computational perspective. We describe an algorithm, and two independent implementations, to compute Bianchi modular forms and their Hecke eigensystems over general imaginary quadratic fields, extending the computations done by the first author and several of his students \cite{cremonathesis,cremona_hyptess,whitley,bygott,lingham,aranes,cremona_aranes} over fields of class number $1$, $2$, and $3$. We give an overview of the results obtained from these implementations, together with details for the field $\QQ(\sqrt{-17})$, whose class group is cyclic of order $4$.

Computations of Bianchi modular forms date back to the 1980s in the work of Grunewald, Mennicke, and others \cite{grunewaldhellingmennicke,elstrodtgrunewaldmennicke}.  They computed Bianchi modular forms at prime levels for $K=\QQ(\sqrt{-d})$ where $d=1,2,3$. In the first author's thesis and 1984 paper \cite{cremonathesis, cremona_hyptess},  these computations were extended to all five Euclidean fields and arbitrary levels, using modular symbols. In the years that followed, further extensions were made to the modular symbol method by several of his students:  Whitley's thesis \cite{whitley} covered the remaining four fields of class number $1$; Bygott's thesis \cite{bygott}  developed techniques for computing Bianchi modular forms over imaginary quadratic fields whose class group is an elementary abelian $2$-group, with explicit examples for $K=\QQ(\sqrt{-5})$ with class number 2; Lingham's thesis \cite{lingham} considered the odd class number case, and computed explicit examples for the fields $\QQ(\sqrt{-23})$ and $\QQ(\sqrt{-31})$ of class number $3$; Aranes's thesis \cite{aranes} made progress towards extending modular symbol methods to imaginary quadratic fields with arbitrary class group.  More recently, the first author has developed a complete implementation that works on arbitrary imaginary quadratic fields in his {\tt C++} package {\tt bianchi-progs}, available in the \Github\ repository \cite{Cremona_bianchi_progs}, and the resulting data is tabulated in the L-functions and modular forms database (\lmfdb) \cite{lmfdb}.

The aforementioned work exploits the connection between Bianchi modular forms and the homology of certain quotients of the hyperbolic $3$-space $\hy_3$. To compute homology, we require a tessellation of the hyperbolic 3-space $\hy_3$, on which the Bianchi group $\GL{\OO}$ and its congruence subgroups act.  One way to obtain such a tessellation, used by the first author and his students, utilizes an algorithm coming from the work of Swan \cite{swan}. The tessellation data from Swan's algorithm can also be used to give a ``pseudo-Euclidean algorithm'' for all imaginary quadratic fields $K$, first introduced in \cite{whitley}, which plays the part of the usual continued fraction algorithm as used in \cite{cremona_hyptess}.   An alternative approach was developed by the third author and uses the work of Ash \cite{ash-def} and Koecher \cite{koecher-cones} coming from the theory of perfect Hermitian forms. 
We note that while both methods involve first finding a tessellation of $\hy_3$ by ideal polyhedra (a precomputation which only needs to be done once for each field), the two methods produce (in general) different tessellations.
For Hecke operator computations, the second method uses the reduction theory introduced by Gunnells \cite{gunnells-qrank1} instead of the pseudo-Euclidean algorithm.

Using the two methods just described, we have two independent implementations of algorithms to compute spaces of Bianchi modular forms, and the action of the Hecke algebra on them, for arbitrary imaginary quadratic fields;  in both cases there is a precomputation stage which only needs to be carried out once for each field (computing the tessellation data), then the rational homology $H_1(\hy_3/\Gamma_0({\nn}), \QQ)$ is computed, and finally the Hecke action on this space and its decomposition into Hecke eigenspaces.   The first implementation is in the {\tt C++} package {\tt bianchi-progs} \cite{Cremona_bianchi_progs} mentioned above;  the second method  has been implemented in \Magma\ \cite{MAGMA} by the third author (see \cite{yasak_hyptess}) and extended by the second author in her thesis \cite{kalanithesis}.   Details of the Bianchi newforms spaces, and of individual Bianchi newforms (currently only those with rational Hecke eigenvalues) may be found in the \lmfdb.

The Bianchi modular forms we compute using the $1$-homology of quotients of hyperbolic $3$-space $\hy_3$ include forms with unramified character, as well as those of trivial character, which have been of greatest interest to date.  They have a weight representation, which is analogous to weight $2$ for classical modular forms, given by a $3$-dimensional irreducible representation of $\SU(2,\CC)$, and these Bianchi modular forms are vector-valued functions, taking values in $\CC^3$, and corresponding to harmonic differential $1$-forms on quotients of $\hy_3$. (Analogously, the weight of a classical modular form is an integer $k$ which determines a $1$-dimensional representation of the group $\SO(2,\RR)$, and holomorphic cusp forms of weight $2$ correspond to holomorphic differential $1$-forms on quotients of the upper half-plane $\hy_2$.)  Eigenforms with trivial character and rational eigenvalues are related to elliptic curves (or, in some cases, Abelian varieties of higher dimension); hence, for each such eigenform we compute, we have also attempted to find an elliptic curve which matches it (in the sense of having the same Galois representation and the same $L$-function).  Further work is then required to prove that the Bianchi newform and the elliptic curve really do have isomorphic Galois representations.  In one of the examples we give for the field $K=\QQ(\sqrt{-17})$ (see Example~\ref{ex:modularity-example}), we prove the modularity of an elliptic curve over $K$, the curve with \lmfdb\ label~\lmfdbecnf{2.0.68.1}{7.2}{a}{2}. This example complements the work of by Dieulefait, Guerberoff, and  Pacetti where they prove the modularity of certain explicit elliptic curves over imaginary quadratic fields $\QQ(\sqrt{-d})$ for $d=23, 31$ which have class number $3$ \cite{modularityexamplesDGP}. We note that the recent paper of Newton and Caraiani \cite{caraiani2025modularityellipticcurvesimaginary} establishes the modularity of all elliptic curves over a subset of imaginary quadratic fields, namely, those $K$ for which the Mordell-Weil rank of the elliptic curve $X_0(15)(K)$ is zero.  However, while this general result does include both $\QQ(\sqrt{-23})$ and $\QQ(\sqrt{-31})$, it does not include $K=\QQ(\sqrt{-17})$, since $X_0(15)(\QQ(\sqrt{-17}))$ has rank $2$.

Where possible, we have carried out all computations reported on here independently using both implementations, to check for consistency.  This consistency check has been carried out for (at least) the following:
\begin{itemize}
\item  homology dimensions, and integral homology structure, at level~$\nn=(1)$ for all fields $K$ with $|\disc(K)|\le2100$;
\item  homology dimensions for levels $\nn$ with $\n(\nn)\le100$, for fields $K$ with $|\disc(K)|\le100$;
\item homology and cusp form space dimensions, and Hecke eigenvalues and eigensystems for $\n(\nn)\le200$ for $K=\QQ(\sqrt{-17})$.
\end{itemize}
For example, for the field $K=\QQ(\sqrt{-17})$, the dimensions of the spaces of Bianchi cuspforms at levels ${\nn}$ with $N({\nn})\le1000$ were computed using the {\tt C++} implementation, and may be be found in the \lmfdb\ (the table may be seen at \url{https://www.lmfdb.org/ModularForm/GL2/ImaginaryQuadratic/gl2dims/2.0.68.1}).  The dimensions for $N({\nn})\le200$ were also computed using the \Magma\ implementation, and found to agree.

Both our implementations include code for computing all (principal and non-principal) Hecke operators on homology, and hence obtaining full Hecke eigensystems,  using results from the first author's paper \cite{cremona-preprint}. 

The structure of this paper is as follows. In Section~\ref{sec:BMF}, we briefly recall the notions of Bianchi modular forms and adelic Bianchi modular forms. In Section~\ref{sec:eigensystems}, we discuss Hecke eigensystems, homological eigensystems, and how to recover a complete Hecke eigensystem from a homological eigensystem. In Section~\ref{sec:resultsandtables}, we provide dimension tables and explicit detailed examples from the implementation of the algorithm for the imaginary quadratic field $K=\QQ(\sqrt{-17})$. In a final short section we give a brief survey of the data we have so far computed for other imaginary quadratic fields.

The authors wish to thank Paul~Gunnells and John~Voight for many helpful conversations, and Frank~Calegari for explaining to the first author the reason for multiple self-twists being impossible (see Remark~\ref{rmk:at-most-one-self-twist} below). The third author was supported in part by a grant from the Simons Foundation (848154, DY). 

\section{Bianchi Modular Forms}
\label{sec:BMF}

As a general reference for Bianchi groups and Bianchi modular forms, we refer to {\c Seng\"un}'s survey article \cite{sengun2014arithmetic}, which has a comprehensive bibliography.  Bygott's thesis \cite{bygott} also contains many details about the different ways to define Bianchi modular forms (as vector-valued functions on $\hy_3$, as adelic functions, and as functions on ``modular points'') and how these are related.  Some of the facts we use about Bianchi modular forms are special cases of a much more general theory of automorphic forms over general global fields, as developed in Weil's book \cite{Weil} and Miyake's paper \cite{Miyake}.

\subsection{Bianchi modular forms and adelic Bianchi modular forms}
Let $K$ be an arbitrary imaginary quadratic field, with class group $\CL$ and class number $\CLno=|\CL|$.
The Bianchi modular forms we are mainly concerned with here are cuspidal forms of ``weight\footnote{As remarked in the introduction, the weight is actually determined by a $3$-dimension representation of $\SU(2,\CC)$.  We refer to this as ``weight $2$'' for simplicity, as they are the closest analogue of classical forms of weight $2$, and include base changes of classical weight $2$ forms  from $\QQ$ to $K$.} $2$'' and level $\Gamma_0(\nn)$, with trivial character.  Here, $\nn$ is an integral ideal of $K$, and $\Gamma_0(\nn)$ is the congruence subgroup
\[
\Gamma_0({\nn}) = \set*{ \mat{a&b\\c&d} \in \go \given c \in {\nn} }.
\]
One initially defines Bianchi modular forms as vector-valued functions 
\[F=(F_0,F_1,F_2)\colon \hy_3 \to \CC^3\] satisfying certain analytic conditions, which amount to saying that the associated differential $1$-form $F_0(z,t)\frac{-dz}{t}+F_1(z,t)\frac{dt}{t}+F_2(z,t)\frac{d\overline{z}}{t}$ is harmonic and invariant under the action of $\Gamma_0(\nn)$ on $\hy_3$, so defines a harmonic differential on the quotient real analytic space $X_0(\nn) := \hy_3^*/\Gamma_0(\nn)$; here $\hy_3^*=\hy_3\cup K\cup\{\infty\}$ is the completion of $\hy_3$ obtained by including the $K$-rational cusps $\PP^1(K)=K\cup\{\infty\}$.   Such functions have Fourier expansions in terms of $K$-Bessel functions, in the form of series expansions which are sums of terms indexed by non-zero elements of the ring of integers $\OO$. The coefficients for associate elements of $\OO$ are equal\footnote{This would not be the case if we worked with subgroups of $\SL(2,\OO)$: see \cite{cremonathesis}.}, so this is expansion can also be written as a sum over principal ideals.  However, in order to obtain a space on which the full Hecke algebra $\TT$ acts, and for compatibility with the definition of ``adelic Bianchi modular forms'' which are functions on the adelic space $\GL{\AA_K}$, it is necessary to consider $\CLno$-tuples of such functions, or equivalently, functions on the union $\widetilde{\hy_3} = \cup_{c\in\CL}\hy_3^{(c)}$ of $\CLno$ copies of $\hy_3$, where the function on each copy $\hy_3^{(c)}$ is invariant under a twisted form $\Gamma_0^{(c)}(\nn)$ of $\Gamma_0(\nn)$.  This group is a subgroup of $\GL{K}$ which depends on a choice of an ideal $\pp$ in the class $c$:
\[
  \Gamma_0^{\pp}(\nn) = \set*{\mat{a&b\\c&d}\given a,d\in \OO;
  b\in\pp^{-1}; c\in\nn\pp; ad-bc\in \OO^\times};
\]
it is a form of adelic conjugate of $\Gamma_0(\nn)$.  We can write such a function as $F=\sum_{c\in\CL}F^{(c)}$, where each component $F^{(c)}\colon \hy_3^{(c)}\to\CC^3$ transforms under the action of $\Gamma_0(\nn)^{(c)}$.  These functions determine harmonic differentials on the larger space 
\[
\widetilde{X_0(\nn)} = \bigcup_{c\in\CL}X_0(\nn)^{(c)} = \bigcup_{c\in\CL} \hy_3^{*(c)}/\Gamma_0(\nn)^{(c)};
\]
they also have Fourier expansions, which are sums over all integral ideals $\aa$ of $\OO$, with coefficients $a_F(\aa)$.  The terms indexed by ideals in each ideal class $c$ come from the expansion of the component $F^{(c)}$.  Following {\c Seng\"un}, we will call these functions $F\colon \widetilde{\hy_3}\to\CC^3$ \emph{adelic Bianchi modular forms}, keeping the simpler term \emph{Bianchi modular forms} for single functions $\hy_3\to\CC$.   We write $\hy_3$ for the principal component $\hy_3^{(1)}$ and $X_0(\nn)=X_0(\nn)^{(1)}$.

Denote the space of (cuspidal, weight $2$) Bianchi modular forms of level $\Gamma_0(\nn)$ and character $\chi$ by $S(\nn,\chi)$, let $S(\nn)=S(\nn,1)$, the subspace of those with trivial character, and $S(\nn)^{\ur}$ the sum of the spaces $S(\nn,\chi)$ over all unramified characters $\chi$.  We identify the group of all unramified characters with the group of characters of $\CL$; they form an abelian group dual to, and hence (non-canonically) isomorphic to, $\CL$.  The class group acts on $S(\nn)^\ur$ via the Hecke operators $T_{\aa,\aa}$, whose action on $S(\nn)^\ur$ only depends\footnote{This is because $\chi$ is unramified, and we are restricting to forms of weight $2$.} on the ideal class $[\aa]$; we have
\[
S(\nn,\chi) = \{F\in S(\nn)^\ur \mid T_{\aa,\aa}(F) = \chi(\aa) F\ \text{for all integral $\aa$ coprime to $\nn$}\};
\]
here, and below, we write $\chi(\aa) = \chi(c)$ where $c=[\aa]\in\CL$ is the class of the ideal $\aa$. In particular,
\[
S(\nn) = \{F\in S(\nn)^\ur \mid T_{\aa,\aa}(F) = F\ \text{for all integral $\aa$ coprime to $\nn$}\}.
\]
Twisting by the unramified character $\psi$ maps $F\in S(\nn,\chi)$ to $F\otimes\psi\in S(\nn,\chi\psi^2)$.  If $F\in S(\nn,\chi)^\new$ (see below for the definition of the new subspace) then writing $F = \sum_{c\in\CL}F^{(c)}$, we have $F\otimes\psi=\sum_{c\in\CL}\psi(c)F^{(c)}$.

It is possible for an adelic Bianchi modular form to be its own twist by a nontrivial unramified character $\psi$;  this requires $\psi^2$ to be trivial, so can only occur when $\CLno$ is even, as otherwise no unramified quadratic characters exist.  Such a form is said to admit a \emph{self-twist} by the character $\psi$.  We will see examples of forms which admit self-twist below; their existence presents certain computational difficulties which we must overcome.  
Note that when $F$ admits self-twist by $\psi$, half of its component functions (namely, the $F^{(c)}$ for which $\psi(c)=-1$) are identically zero.  If $F$ is an eigenform for the Hecke operator $T_\aa$ with eigenvalue $\alpha(\aa)$ then $F\otimes\psi$ is also an eigenform with eigenvalue $\psi(\aa)\alpha(\aa)$.  Hence, if $F$ has self-twist by $\psi$, we have $T_\aa(F)=0$ for all $\aa$ with $\psi(\aa)=-1$.

The full Hecke algebra $\TT$ is generated by operators $T_\aa$ (for all integral ideals $\aa$) and $T_{\aa,\aa}$ (for integral ideals $\aa$ coprime to the level).  The algebra has a grading by the class group, where the class of $T_\aa$ is $[\aa]$ and that of $T_{\aa,\aa}$ is $[\aa^2]$.  Following \cite{cremona-preprint}, we denote by $\TT_{c}$ the subset of operators of class $c$, for each $c\in\CL$, with $\TT_{1}$ the subalgebra of principal Hecke operators; the grading means that $\TT_{c}\TT_{c'}\subseteq\TT_{cc'}$.  If $F=\sum_cF^{(c)}\in S(\nn)^\ur$ and $T(F)=G$ for some Hecke operator in class $[T]$, we have $T(F^{(c)}) = G^{([T]c)}$ for each $c\in\CL$.  It follows that the principal subalgebra $\TT_{1}$ acts on the space of principal components $\{ F^{(1)} \mid F \in S(\nn)^\ur\}$.

An adelic Bianchi modular form which is an eigenvector for all Hecke operators will be called an \emph{eigenform}.
As with classical modular forms, the Fourier coefficients of an eigenform $F$ are determined by its Hecke eigenvalues.  The new subspace $S(\nn)^{\ur,\new}$ (defined below), has a basis consisting of eigenforms, and these have nonzero unit coefficient $a_F(1)$; we normalise them to have unit coefficient equal to $1$, and call normalised new eigenforms \emph{newforms}.   For a newform $F$, the eigenvalue of $T_\pp$ is equal to $a_F(\pp)$.  The coefficients are multiplicative, and in general satisfy the same multiplicative relations as the Hecke operators themselves, so that the $\aa$-coefficient is the eigenvalue of $T_\aa$ for all integral ideals $\aa$ coprime to the level $\nn$.  We refer to Miyake \cite{Miyake} for these properties of adelic modular forms.

\begin{remark}\label{rmk:at-most-one-self-twist}
    An newform can admit at most one self-twist, even when the class group admits more than one quadratic character.  This follows from the existence of Galois representations attached to such forms.  (It is known that these Bianchi modular forms have associated
    $\ell$-adic Galois representations $\rho_{F,\ell}$.  These were first
    constructed by Taylor \textit{et al.}~in \cite{TaylorI},
    \cite{TaylorII} with subsequent results by Berger and Harcos
    in~\cite{Berger}.)   The reason (as explained to the first author by Frank Calegari) is essentially that if a newform had two distinct self-twists, then the associated projective Galois representation would have finite image, contradicting known facts about these representations. 
\end{remark}

A fundamental result of Kur\v{c}anov \cite[Theorem~1]{kurcanov}  is that the complex cohomology (and hence also the complex homology) of the quotient space $\widetilde{X_0(\nn)}$ is isomorphic to $S(\nn)^\ur$, the space of all cuspidal Bianchi modular forms on $\Gamma_0(\nn)$ with unramified character. More precisely, integration of harmonic differentials along closed paths induces a perfect dual pairing between these two complex vector spaces, which respects the Hecke action on both sides.  Using the Hecke operators $T_{\aa,\aa}$, we can cut out a subspace of this homology space whose complexification is isomorphic to the space of Bianchi modular forms with trivial, or any other unramified, character.

For computational purposes, it is more convenient to compute only the homology of the principal part of the quotient space, that is,  $H_1(X_0(\nn), \QQ) = H_1(\hy_3^*/\Gamma_0(\nn), \QQ)$.  On this smaller principal space, the full Hecke algebra $\TT$ does not act, as an operator $T\in\TT_{c}$ maps $H_1(X_0(\nn)^{(c')}, \QQ)$ to $H_1(X_0(\nn)^{(c'c^{-1})}, \QQ)$.  Only the principal  subalgebra $\TT_{1}$ acts on the principal homology $H_1(X_0(\nn), \QQ)$ itself.  The principal operators include $T_{\aa}$ for principal ideals $\aa$ coprime to $\nn$ and $T_{\aa,\aa}$ for $\aa$ coprime to $\nn$ with $\aa^2$ principal, and more generally $T_{\aa,\aa}T_{\bb}$ where $\aa$ and $\bb$ are coprime to $\nn$ and $\aa^2\bb$ is principal.  For explicit formulas for these, we refer to the first author's paper \cite[\S 4]{cremona-preprint}, which builds on special cases developed by Aranes \cite{aranes}, Bygott \cite{bygott}, Lingham \cite{lingham}, and Whitley \cite{whitley}.   The approach used there is to view Bianchi modular forms as functions on \emph{modular points}, generalising the treatment in the classical case (over $\QQ$) of Serre for level $1$ in \cite[VIII~\S1]{serre-book} and for arbitrary levels by  Koblitz in \cite[III~\S5]{KoblitzECMF} and Lang \cite[VII]{LangModularForms}.   This extends the work of Bygott in his thesis \cite{bygott}, where  a theory of modular points for general number fields was first developed.

\subsubsection{Atkin-Lehner operators}
As well as the Hecke operators $T_{\aa}$ and $T_{\aa,\aa}$, we also have Atkin-Lehner operators $W_\qq$ operating on $S(\nn)^\ur$ for each $\qq\mid\nn$ such that $\qq$ and $\qq^{-1}\nn$ are coprime.   Such divisors of the level are called \emph{exact divisors}, and we write $\qq\mid\mid\nn$. In general, we have $W_{\qq}^2=T_{\qq,\qq}$ (by \cite[\S3.2.1]{cremona-preprint}), so their restrictions to the trivial subspace $S(\nn)$ they are (as in the classical case)  involutions,  forming an elementary abelian $2$-group of rank equal to the number of prime divisors of $\nn$.  This group is generated by the $W_\qq$ for $\qq$ an exact prime power divisor of $\nn$.  Also restricting to $S(\nn)$, Atkin-Lehner operators commute with Hecke operators $T_{\aa,\aa}$ and $T_\aa$ for $\aa$ coprime to $\nn$  by \cite[Prop.3.5]{cremona-preprint}, so take eigenforms in $S(\nn)$  to eigenforms; in particular, every eigenform for the Hecke algebra $\TT$ with trivial character is also an eigenform for each Atkin-Lehner operator $W_\qq$ with eigenvalue $\pm1$.  Twisting takes eigenforms to eigenforms;  as with the $T_\aa$ operators, twisting an eigenform with trivial character by an unramified quadratic character $\psi$ multiplies the $W_\qq$-eigenvalue by $\psi(\qq)$.

We will not consider here the action of Atkin-Lehner operators on forms with nontrivial character, where in general (as in the case of classical modular forms, see \cite{AtkinLi}) they map eigenforms to different eigenforms, and one is led to defining their ``pseudo-eigenvalues'' instead of eigenvalues.

\begin{proposition}
  \label{prop:AL-restricts-self-twists}
      For each integral ideal $\nn$, let $C(\nn)$ be the (possibly empty) set of quadratic characters $\psi$ of $\CL$ such that $\psi(\qq)=+1$ for all proper divisors $\qq\mid\mid\nn$.  If an eigenform $F\in S(\nn)^{\ur}$ admits a non-trivial self-twist, then the self-twist character must belong to $C(\nn)$.
\end{proposition}
\begin{proof}
    Suppose that $F$ is an eigenform with self-twist by $\psi$. Let $W_\qq(F)=\varepsilon_F(\qq)F$.   Since $W_\qq^2=T_{\qq,\qq}$ is invertible,  $\varepsilon(\qq)\not=0$, and since $F=F\otimes\psi$, we have $\varepsilon_F(\qq)=\psi(\qq)\varepsilon_F(\qq)$, so $\psi(\qq)=+1$.
\end{proof}
In our algorithm, for each level $\nn$ we determine the set $C(\nn)$ of ``eligible'' self-twist characters for eigenforms at level $\nn$, and this can help to eliminate the (relatively rare) possibility of the existence of self-twist eigenforms at that level.

As with Hecke operators, the Atkin-Lehner operator $W_\qq$ only acts on principal homology when $\qq$ is principal.  In order to compute all Atkin-Lehner eigenvalues, we will also use the principal operators $T_{\aa,\aa}W_\qq$  (when $\aa^2\qq$ is principal), and $T_\aa W_\qq$ (when $\aa\qq$ is principal), again using explicit formulas from \cite{cremona-preprint}.

\begin{example}
    \label{ex:ur-twists-17}
    Let $K=\QQ(\sqrt{-17})$, with $\CL$ cyclic of order $4$.  Let $c$ be an ideal class generating $\CL$ (for example, one can take $c$ to be the class of the prime ideal $(3,1+\sqrt{-17})$).   The character group is also cyclic, generated by the character $\chi_1$ of order $4$ defined by $\chi_1(c)=\sqrt{-1}$.  There is exactly one unramified quadratic character $\chi_2=\chi_1^2$, associated to the unramified quadratic extension $K(\sqrt{-1})$, the genus field of $K$.  The rational primes $p$ which split as $(p)=\pp\overline{\pp}$ in $K$ are those with $\leg{-17}{p}=+1$,  and    for these primes we have $\chi_2(\pp)=+1$ if and only if also $p\equiv1\pmod{4}$.  For both the ramified primes $\pp$ (above $2$ and $17$) we have $\chi_2(\pp)=+1$.  Hence, the set of primes $\pp$ with $\chi_2(\pp)=-1$ (equivalently, those whose ideal classes have order $4$) consists of those above rational primes $p$ such that $p\equiv3\pmod4$ and $p\equiv\pm3,\pm5,\pm6,\pm7\pmod{17}$.  Over this field, therefore,  forms with non-trivial self-twist can only exist at level $\nn$ if all prime factors $\pp$ of $\nn$ with norm in this set divide $\nn$ to an even power.  The first example of a self-twist newform with trivial character occurs at level $\nn=(8)$.  It is the base-change of the classical cusp form with \lmfdb\ label \lmfdbcmf{32}{2}{a}{a} (the unique newform of weight $2$, trivial character, and level $32$), associated to the elliptic curve \lmfdbec{32}{a}{1}, which has complex multiplication by $\ZZ[\sqrt{-1}]$. 
\end{example}

\subsubsection{Newforms and oldforms}
\label{sec:newform-oldforms}
Let $\mm$ and $\nn$ be two levels, with $\mm\mid\nn$, $\mm\not=\nn$.  For each ideal divisor $\dd$ of $\mm^{-1}\nn$, there is an operator $A_\dd \colon S(\mm)^\ur\to S(\nn)^\ur$, which commutes with the Hecke action of Hecke operators $T_\aa$ (for $\aa$ coprime to $\nn$) and $T_{\aa,\aa}$ (see \cite[\S4.4]{cremona-preprint}).   These are the analogues for adelic Bianchi modular forms of the classical maps taking a function $f(\tau)$ on the upper half-plane to $f(d\tau)$, where $d$ is a positive divisor of the quotient of two levels.  These maps take eigenforms to eigenforms, with the same eigenvalues away from the level.  The images of an eigenform $G\in S(\mm)^\ur$ for distinct divisors $\dd$ are linearly independent, and hence span a subspace of $S(\nn)^\ur$ called the \emph{oldclass} associated to $G$, whose dimension is therefore equal to $\sigma_0(\mm^{-1}\nn)$, where $\sigma_0(\aa)$ denotes the number of ideal divisors of an integral ideal $\aa$.  The \emph{oldspace} at level $\nn$ is  the sum of all the oldclasses as $\mm$ ranges over all  levels $\mm\mid\nn$, $\mm\not=\nn$.  The \emph{newspace} $S(\nn)^{\ur,\new}$ is the orthogonal complement of the oldspace with respect to an inner product on $S(\nn)^\ur$ which generalises the classical Peterssen inner product (for a definition of this inner product for adelic Bianchi modular forms, see Miyake \cite{Miyake}). Recall that we define a newform to be an eigenform in the new subspace.  We summarise here some key facts from Miyake's paper, noting that Miyake deals with forms with arbitrary character (of conductor not dividing the level), and defines an eigenform at level $\nn$ to be a form which is an eigenvector for all $T_\pp$ and also for all the adjoint operators $T_\pp^*$;  for the subspace $S(\nn)^\ur$, this is equivalent to being an eigenvector for all $T_\pp$ and all $T_{\pp,\pp}$ for $\pp\nmid\nn$, since $T_{\aa}^* = T_{\aa,\aa}^{-1}T_{\aa}$ (see \cite[\S3.4.2]{cremona-preprint}).
\goodbreak
\begin{theorem}
\label{thm:miyake-newforms}\mbox{}
    \begin{enumerate}
    \item The Hecke algebra $\TT$ is commutative; its restriction to $S(\nn)^{\ur,\new}$ is semisimple. Each Hecke operator is diagonalisable. Hence $S(\nn)^{\ur}$ has a basis of eigenforms, and $S(\nn)^{\ur,new}$ has a basis of newforms (see \cite[Theorem 1]{Miyake}).
    \item For every eigenform $F\in S(\nn)^\ur$,  there is a newform $G\in S(\mm)^{\ur,\new}$ for some level $\mm\mid\nn$ with the same Hecke eigenvalues $T_{\pp}$ for $\pp\nmid\nn$ as $F$, so that either $F$ is itself a newform, or, if $\mm\not=\nn$, $F$ is in the oldclass associated to $G$ (see \cite[Corollary to Theorem 1]{Miyake}).
    \item If newforms $F_j\in S(\nn_j)^{\ur,\new}$ for $j=1,2$ have the same eigenvalues for $T_\pp$ for almost all primes $\pp$, then $\nn_1=\nn_2$ and $F_1=F_2$  (see \cite[Theorem 2]{Miyake}).
\end{enumerate}
\end{theorem}
Hence the multiplicity at level $\nn$ of the eigensystem associated to an eigenform $G\in S(\mm)^{\ur,\new}$ is $\sigma_0(\mm^{-1}\nn)$.
This leads to an expression for the dimensions of the newspaces $S(\nn)^{\ur,\new}$ in terms of the dimensions of the full spaces $S(\nn)^\ur$:
\begin{equation}
\label{eqn:newspace-dimension-formula}
\dim S(\nn)^{\ur,\new} = \dim S(\nn)^{\ur} - \sum_{\mm\mid\nn, \mm\not=\nn}\sigma_0(\mm^{-1}\nn)\dim S(\mm)^{\ur,\new}.
\end{equation}

When we consider the dimensions of the projections from $S(\nn)^\ur$ to the principal component, the corresponding  formula requires some adjustment to take into account eigenforms with self-twist; this phenomenon was first remarked in Bygott's thesis \cite{bygott}.  The reason is as follows.  The maps $A_\dd$ have class $[\dd]$, in the sense that if $G\in S(\mm)^\ur$ and $A_\dd(G)=F\in S(\nn)^\ur$, (where $\mm\mid\nn$ and $\dd\mid\mm^{-1}\nn$), writing $F=\sum_cF^{(c)}$ and $G=\sum_cG^{(c)}$ we have $A_\dd(G^{(c)}) = F^{([\dd]c)}$.  In most cases, a new eigenform $G$ will have all components $G^{(c)}$ non-zero (see Proposition~\ref{prop:zero-components} below).  If this were always the case, then a formula similar to (\ref{eqn:newspace-dimension-formula}) would still hold for the dimensions of the principal components.  However, when $G$ is a new eigenform with self-twist by the unramified quadratic character $\psi$, its components $G^{(c)}$ satisfy $G^{(c)} \not=0$ if and only if $\psi(c)=+1$;  in particular, the principal component $G^{(1)}$ is nonzero.   Now, when applying the operator $A_\dd$, we see (Corollary~\ref{cor:zero-components}) that the principal component of $A_\dd(G)$ is nonzero if and only if $\psi(\dd)=+1$.  Hence, the contribution of $G$ to the principal oldspace dimension at level $\nn$ is not $\sigma_0(\mm^{-1}\nn)$, the total number of divisors of $\mm^{-1}\nn$, but is instead the possibly smaller number of divisors $\dd$ such that $\psi(\dd)=+1$.  For example, if $\nn=\mm\pp$ with $\pp$ prime, then we expect the dimension of the oldspace at level $\nn$ associated to a new eigenform at level $\mm$ to be $2$, the number of ideal divisors of $\pp$, but if $G$ admits self-twist by $\psi$ and $\psi(\pp)=-1$, then the multiplicity is only $1$. This observation can be used to prove that a newform admits a self-twist, by computing the dimension of the associated oldspace at a suitably chosen higher level.

In summary:
\begin{proposition}
\label{prop:oldform-multiplicities}
    Let $G\in S(\mm)^{\ur,\new}$ be a newform, let $\nn$ be a multiple of $\mm$, and $V\subseteq S(\nn)^\ur$ be the oldclass associated to $G$ at level $\nn$, and let $V^{(1)}$ be the projection of $V$ to the principal component.  Then
    \begin{enumerate}
        \item $\dim V = \sigma_0(\mm^{-1}\nn) = \sum_{\dd\mid\mm^{-1}\nn}1$.
        \item
        \begin{enumerate}
            \item If $G$ admits no self-twist then $\dim V^{(1)} = \dim V$.
             \item If $G$ admits self-twist by the character $\psi$ then 
             \[
             \dim V^{(1)} = \#\{\dd\mid\mm^{-1}\nn \mid \psi(\dd)=+1\} = \sum_{\dd\mid\mm^{-1}\nn}(1+\psi(\dd))/2.
             \]
        \end{enumerate}
    \end{enumerate}
\end{proposition}

This makes it impossible to give a simple formula such as (\ref{eqn:newspace-dimension-formula}) for the new dimensions of the principal components in general, when $\CLno$ is even.   It is only possible in certain special cases, for example when $\psi(\pp)=+1$ for all unramified quadratic characters $\psi$ and all primes $\pp\mid\nn$. Otherwise, we can only compute the dimension of the principal component of $S(\nn)^{\ur,\new}$ when we know not just the dimensions $\dim S(\mm)^{\ur,\new}$ for proper divisors $\mm$ (that is, the total number of newforms at level $\mm$), but also, in each case, the number of newforms which have self-twist by each unramified quadratic character $\psi$.

\bigskip

\section{Eigenforms and Eigensystems}
\label{sec:eigensystems}
Our implementations rely on the fact that knowing only the principal homology, and the action on this of principal Hecke operators, suffices to determine the complete set of Hecke eigensystems for the full Hecke algebra $\TT$ acting on the full space with all $\CLno$ components.  That this is the case (and even in a more general setting than just Bianchi modular forms) is a main result of \cite{cremona-preprint}:
\begin{theorem}[Theorem 4.1 of \cite{cremona-preprint}]
\label{thm:preprint-3.7}
   Two eigenvalue systems $\lambda_1$ and $\lambda_2$ of the same level
  $\nn$ have the same restriction to the subalgebra $\TT_{1}$ of principal
  operators if and only if $\lambda_2=\lambda_1\otimes\psi$ for some unramified 
  character~$\psi$ (that is, for some character of the ideal class group~$\CL$).
  The characters~$\chi_1$, $\chi_2$ of $\lambda_1$, $\lambda_2$
  satisfy $\chi_2=\chi_1\psi^2$, and hence $\lambda_1$ and $\lambda_2$
  have the same character if any only if $\psi$ is quadratic.  
\end{theorem}

Hence, knowledge of the eigenvalues of only the principal operators is enough, in principle,  to characterize a full eigensystem, up to unramified twist.   There remains the algorithmic question of actually computing the eigenvalues of non-principal operators (for one, and hence all, of this family of twists) from those of the principal operators.  A solution to this was sketched in \cite[Remark 6]{cremona-preprint}, and we give further details here, as this is an important part of our algorithm.  We express the algorithm in terms of eigensystems, but it can also be seen as a method to recover an eigenform (up to twist) from its principal component.
 
This process is very much simpler over fields with odd class number, and we describe this first, before describing the general situation, where the genus group $\CL/\CL^2$ plays a role.  This approach for fields of odd class number was first used in Lingham's thesis \cite{lingham}.

\subsubsection*{The odd class number case}
Suppose that the class number $\CLno$ is odd.   Then every element of the class group is a square, every unramified character is a square, and there are no unramified quadratic characters.  It follows that each space $S(\nn,\chi)$ is isomorphic to $S(\nn)$;  in particular, $\dim S(\nn)^\ur = \CLno\dim S(\nn)$.  Taking $\psi$ to be a character with $\psi^2=\chi$, twisting by $\psi$ gives an isomorphism $S(\nn)\to S(\nn,\chi)$.  In other words, every eigenform in $S(\nn)^\ur$ may be obtained from an eigenform with trivial character by twisting.   Secondly, when the ideal class $c$ is a square, then the twisted group $\Gamma_0(\nn)^{(c)}$ is actually conjugate by a matrix $M\in\GL{K}$ to the principal group $\Gamma_0(\nn)$. (Essentially, $M$ is a matrix inducing the isomorphism $\aa^2\oplus\OO\cong\aa\oplus\aa$ for a suitable ideal $\aa$ with $c=[\aa]^2$.)  It follows that $M$ induces an isomorphism $H_1(X_0(\nn)^{(c)},\QQ) \cong H_1(X_0(\nn), \QQ)$.  Hence all these spaces have the same dimension (as $\QQ$-vector spaces), and, by Kur\v{c}anov's theorem, this dimension is equal to $\dim_{\CC} S(\nn)$.   Moreover, we can compute the action of every $T_{\bb}$ on $H_1(X_0(\nn), \QQ)$ (for $\bb$ coprime to $\nn$) by computing that of the principal operator $T_{\aa,\aa}T_\bb$, where $\aa$ is an ideal coprime to $\nn$ such that $\aa^2\bb$ is principal, using the formula for such operators given in \cite{cremona-preprint}.  Then the eigenvalues of $T_\bb$ on $S(\nn)$ may be determined, using the fact that the only principal operator of the form $T_{\aa,\aa}$ is the identity.   On  $H_1(X_0(\nn), \QQ)$, therefore, the eigenvalue of $T_\bb$ (on the full eigensystem with trivial character) is the same as that of the principal operator $T_{\aa,\aa}T_\bb$.

\subsection{Hecke eigensystems}

Hecke operators commute and are diagonalizable, and thus are simultaneously diagonalizable. Hence the space $S(\nn)^\ur$ has a basis consisting of eigenforms for all $T\in\TT$.  To each such eigenform we associate a \emph{Hecke eigensystem}.  As already noted, the operator $T_{\aa,\aa}$ only depends on the ideal class $[\aa]$, so $T_{\aa,\aa}$ is the identity when $\aa$ is principal.

\begin{definition}\label{def:Heckesystem}
Let $F\in S(\nn)^{\ur}$ be an eigenform, that is, a simultaneous eigenvector for all operators in the Hecke algebra $\TT$. The function $\lambda_F\colon \TT \to \CC$, defined by
\[
T(F) = \lambda_F(T)F
\]
for each $T \in \TT$, is a ring homomorphism called the \emph{Hecke eigensystem} associated to $F$. 
For each integral ideal $\aa$ coprime to the level $\nn$, we define $\alpha_F(\aa) = \lambda_F(T_{\aa})$ and $\chi_F(\aa)=\lambda_F(T_{\aa,\aa})$, so that (by definition) $F\in S(\nn,\chi_F)$, and $\chi_F$ is constant on ideal classes, so may be identified with a character of the ideal class group.   It will be convenient to set $\chi_F(\aa)=0$ for ideals $\aa$ which are not coprime to the level.

The \emph{Hecke field} of $F$, denoted $k_F$, is the number field generated by the eigenvalues $\lambda(T)$ for $T\in\TT$.  Since the adjoint of $T_{\aa}$ is $T_{\aa,\aa}^{-1}T_{\aa}$, the eigenvalues of principal operators $T_{\aa}$ are all real, and moreover, if $F\in S(\nn)$ then the Hecke field is totally real.

Since $\TT$ is generated by $T_\aa$ and $T_{\aa,\aa}$ for $\aa$ coprime to $\nn$, the eigensystem is uniquely determined by the values $\alpha_F(\aa)$ and $\chi_F(\aa)$; to emphasize this we will also write the eigensystem $\lambda_F$ associated to $F$ as a pair $\lambda_F=(\alpha_F,\chi_F)$.  The Hecke field is generated by $\alpha_F(\aa)$ and $\lambda_F(\aa)$ for $\aa$ coprime to $\nn$; so, in particular, if $F\in S(\nn,\chi)$ then $k_F$ contains the cyclotomic field $\QQ(\chi)$.
\end{definition} 

While the functions $\lambda_F$ and $\chi_F$ are totally multiplicative, the function $\alpha_F$ is only multiplicative.
\begin{theorem}
\label{thm:mult-relations}
The eigensystem $\lambda_F = (\alpha_F,\chi_F)$ associated to a Hecke eigenform $F\in S(\nn)^\ur$  satisfies the following:
\begin{itemize}
\item If $\aa$, $\bb$ are integral ideals, then
\[
\chi_F(\aa\bb) = \chi_F(\aa)\chi_F(\bb).
\]
\item If $\aa$ and $\bb$ are coprime integral ideals, then 
\[
\alpha_F(\aa)\alpha_F(\bb) =   \alpha_F(\aa\bb).
\]
\item If $\pp$ is a prime ideal, then for all $n \geq 1$ we have
\[
\alpha_F(\pp^n)\alpha_F(\pp) = \alpha_F(\pp^{n+1}) + \n(\pp)\alpha_F(\pp^{n-1})\chi_F(\pp).
\]
(Here we have used the convention that $\chi_F(\pp)=0$ if $\pp\mid\nn$.)
\end{itemize}
\end{theorem}

Since for a normalised eigenform $F$ the $\aa$-Fourier coefficient is equal to the eigenvalue $\alpha_F(\aa)$, it follows that these coefficients satisfy the same multiplicative relations as above.

\begin{definition}
    For each character $\psi$ of the class group, and each eigensystem $\lambda:\TT\to\CC$, the \emph{twist} $\lambda\otimes\psi$ is defined by
    \begin{itemize}
        \item $(\lambda\otimes\psi)(T_\aa) = \psi([\aa])\lambda(T_\aa)$;
        \item $(\lambda\otimes\psi)(T_{\aa,\aa}) = \psi^2([\aa])\lambda(T_\aa)$.
    \end{itemize}
\end{definition}

\begin{proposition}
    If $\lambda=(\lambda_F,\chi_F)$ is the eigensystem associated to the newform $F\in S(\nn,\chi_F)$,  the twisted form $F\otimes\psi$ is also an eigenform, with eigensystem  $\lambda\otimes\psi = (\lambda_F\psi, \chi_F\psi^2) = (\lambda_{F\otimes\psi}, \chi_{F\otimes\psi})$.
\end{proposition}
    \begin{proof}
        Write $F=\sum_cF^{(c)}$, so that $F\otimes\psi=\sum_c\psi(c)F^{(c)}$.

        First consider $T=T_\aa$.  This maps $F^{(c)}$ to the $([\aa]^{-1}c)$-component of $T(F)$, which is $\lambda_F(\aa)F^{([\aa]^{-1}c)}$.  Hence
        \begin{align*}
            T(F\otimes\psi) &= \sum_c\psi(c)T(F^{(c)}) \\
            &= \sum_c\psi(c)\lambda_F(\aa)F^{([\aa]^{-1}c)}\\
            &= \lambda_F(\aa)\psi([\aa]) \sum_c\psi([\aa]^{-1}c)F^{([\aa]^{-1}c)}\\
            &= \lambda_F(\aa)\psi([\aa]) \sum_c\psi(c)F^{(c)}\\
            &= \lambda_F(\aa)\psi([\aa]) (F\otimes\psi).
        \end{align*}

        The case of $T=T_{\aa,\aa}$ is similar, using $T_{\aa,\aa}(F^{(c)}) = \chi_F(\aa)F^{([\aa]^{-2}c)}$.
    \end{proof}

We call the set of Hecke eigensystems obtained by twisting a system $(\lambda,\chi)$ by the characters of the class group the \emph{twist orbit} of $(\lambda,\chi)$, and the set of characters in the twist orbit the \emph{character orbit} of the Hecke eigensystem. The character orbit is one complete coset of the subgroup of squares of class group characters, so has size $|\CL^2|$;  the size of the twist orbit is usually equal to $\CLno$, but if the eigensystem admits self-twist by a  quadratic character $\psi$, the size of the twist orbit is only $\CLno/2$.  In this self-twist case, for each $\aa$, we have $\lambda(\aa)\psi(\aa)=\lambda(\aa)$, and hence $\lambda(\aa)=0$ for all ideals $\aa$ in half of the ideal classes, namely those on which $\psi$ takes the value $-1$.  We use this to help identify Hecke eigensystems with self-twist, for if $\alpha_F(\aa)\not=0$ then only the possible self-twist characters $\psi$ are those with $\psi(\aa)=+1$, if any.

The set of ideal classes $c$ for which there exists a Hecke operator $T\in\TT_{c}$ with $\lambda(T)\not=0$ clearly forms a subgroup of $\CL$, which always contains $\CL^2$ since $\lambda(T_{\aa,\aa})=\chi(\aa)\not=0$ for all $\aa$ coprime to $\nn$.  Denote this subgroup by $H_0(\lambda)$, or by $H_0(F)$ when $\lambda=\lambda_F$.  Usually (and always when $\CLno$ is odd), this subgroup is the whole of $\CL$, but if $\lambda$ admits self-twist by the unramified quadratic character $\psi$, we have $H_0(\lambda)=\ker\psi$, a subgroup of index $2$. Recall from Proposition \ref{prop:AL-restricts-self-twists} that for each level $\nn$, the set $C(\nn)$  of eligible self-twist characters may be a proper subset of the set of all unramified quadratic characters. If, for a given eigensystem, there are  ideals $\aa$ (coprime to the level) whose ideal classes generate $\CL/\CL^2$, such that $\alpha(\aa)\not=0$, then we know that the eigensystem admits no self-twist.    

On the other hand, determining rigorously that an eigensystem does admit a self-twist (when $C(\nn)\not=\emptyset$) is more problematical algorithmically, since even if we find that for some $\psi\in C(\nn)$,  $\lambda(T_\pp)=0$ for many primes $\pp$ such that $\psi([\pp])=-1$, this does not---without further work---imply that $\lambda(T_\pp)=0$ for all such primes.  One method to do this is to use oldform multiplicities, using the observation before Proposition \ref{prop:oldform-multiplicities}; however, this involves working at higher levels, so is inconvenient in practice (and we have not implemented this).

\begin{proposition}
\label{prop:zero-components}
    Let $F = \sum_cF^{(c)} \in S(\nn, \chi)^\new$ be a newform. 
    \begin{enumerate}
        \item If $F$ does not admit self-twist, then every component $F^{(c)}$ is nonzero.
        \item If $F$ admits self-twist by the quadratic character $\psi$, then $F^{(c)}=0$ if and only if $\psi(c)=-1$.
    \end{enumerate}
\end{proposition}
\begin{proof}
For each integral ideal $\aa$ coprime to $\nn$ we have $T_\aa(F)=\lambda_F(\aa)F$, and hence $T_\aa(F^{(c)}) = \lambda_F(\aa) F^{([\aa]^{-1}c)}$, so if $\lambda(\aa)\not=0$ then $F^{(c)}=0\implies F^{([\aa]^{-1}c)}=0$.   Hence the set of classes $c$ for which $F^{(c)}=0$ is a union of cosets of the subgroup $H_0(F)$.

If  $F$ does not admit self-twist, then in every ideal class there exists $\aa$ with $\lambda_F(\aa)\not=0$, and hence every component $F^{(c)}$ is nonzero (as otherwise $F=0$).

If  $F$ admits self-twist by the quadratic character $\psi$, then certainly $F^{(c)}=0$ for classes $c\notin\ker\psi$, and the previous argument shows that $F^{(c)}\not=0$ for all classes $c\in\ker\psi$.
\end{proof}

\begin{corollary}
\label{cor:zero-components}
Let $F = \sum_cF^{(c)} \in S(\mm, \chi)^\new$ be a newform at level $\mm$,  let $\nn$ be a level divisible by $\mm$, let $\dd\mid\mm^{-1}$, and let $F_{\dd} = A_{\dd}(F) \in S(\nn, \chi)$ be the oldform at level $\nn$.
\begin{enumerate}
    \item If $F$ has no self-twist, then all components of $F_{\dd}$ are nonzero.
    \item If $F$ as self-twist character $\psi$, then 
    \[
    F_{\dd}^{c}=0 \iff \psi(c)=-\psi(\dd).
    \]
    In particular, the principal component $F^{(1)}$ is nonzero if and only if $\psi(\dd)=1$.
\end{enumerate}
\end{corollary}
\begin{proof}
These follow from $A_{\dd}(F^{(c)}) = F_{\dd}^{[\dd]c}$.    
\end{proof}

\subsection{Homological eigenforms} 

\begin{definition}
We say that a class $f$ in the principal homology $H_1(X_0(\mf{n});\CC)$ is a \emph{homological eigenform} if it is a simultaneous eigenvector for all principal Hecke operators.     For a homological eigenform $f$ we set $\lambda_f(T)$ to be the eigenvalue of the principal operator $T$ on $f$, so $T(f)=\lambda_f(T)f$.  The map $\lambda_f\colon\TT_{1} \to \CC$ is the \emph{homological eigensystem} associated to $f$.

The \emph{principal} (or \emph{homological}) \emph{Hecke field} $k_f$ is the field generated by all the principal eigenvalues $\lambda_T(f)$.  As principal operators are self-adjoint, this is a totally real number field.  Note that $T_{\aa,\aa}$ is principal if and only if $[\aa]\in\CL[2]$, so that the character associated to a homological eigensystem is quadratic.
\end{definition}

Using our explicit description of the finite-dimensional vector space $H_1(X_0(\mf{n});\QQ)$, and the construction of sets of $2\times2$ matrices with entries in $\OO$ defining the action of all principal Hecke operators from \cite[Section 4]{cremona-preprint} we can compute matrices representing the action of each principal Hecke operator $T$ on $H_1(X_0(\mf{n});\QQ)$. In practice, this requires some form of reduction algorithm to express an arbitrary path between cusps as a sum of edges in our tessellation.  In our implementations, we do this using either the reduction theory developed by Gunnells \cite{gunnells-qrank1}, as implemented in the third author's \Magma\ package, or using a ``pseudo-Euclidean algorithm'' generalising the classical continued fraction algorithm, as developed by the first author and his students (first appearing in Whitley's thesis \cite{whitley}), which is implemented in the {\tt C++} package {\tt bianchi-progs} \cite{Cremona_bianchi_progs}.  Having these matrices, which commute, standard linear algebra techniques allow us to find all the principal eigensystems.

We discard old eigensystems simply by recognising them from previous computations at lower levels.  We can even at this point determine, for each new homological eigensystem $g$ at a lower level $\mm\mid\nn$, the multiplicity of the corresponding eigensystem at level $\nn$, as by Proposition~\ref{prop:oldform-multiplicities} and Corollary~\ref{cor:zero-components} this is $\sigma_0(\mm^{-1}\nn)$, except when $g$ admits self-twist by $\psi$, in which case it is $\frac{1}{2}\sum_{\dd\mid\mm^{-1}\nn}(1+\psi(\dd))$.
Assuming that we carry out our computations systematically in order of level norm, we will already know all the newforms at lower levels, and then may use the formula of Proposition \ref{prop:oldform-multiplicities} to compute the old and new dimensions of $H_1(X_0(\mf{n});\QQ)$. Of course, if the new dimension is zero, then we need carry out no further computations at this level.

Note that it is possible for oldclasses with self-twist to have dimension $1$, so we cannot assume that every $1$-dimensional homological eigenform is new.

\begin{definition}\label{def:matching}
A homological eigenform $f$ \emph{matches} a Hecke eigensystem $\lambda=(\alpha,\chi)$ if 
\begin{enumerate}
    \item $\alpha(\aa) = \lambda_f(T_\aa)$ for each principal ideal $\aa$; 
    \item $\chi(\bb)=\lambda_f(T_{\bb,\bb})$ for each ideal $\bb$ coprime to the level $\mf{n}$ in $\CL[2]$ (that is, such that $\bb^2$ is principal). Note that $\chi(\bb)=\pm1$ for all such $\bb$.
\end{enumerate}
\end{definition}

\subsection{Recovering a full eigensystem from a homological eigensystem}

We now describe a procedure to compute one complete Hecke eigensystem $\lambda=(\alpha,\chi)$ which matches a given homological eigensystem $\lambda_f$ at level $\nn$. The procedure has the following form.
\begin{itemize}
    \item [Input:] For any principal operator $T$, the eigenvalue $\lambda_f(T)$.
    \item [Output:]
    \begin{enumerate}
        \item For any ideal $\aa$ coprime to $\nn$,  the value $\chi(\aa)$.
        \item For any prime $\pp\nmid\nn$,  the value $\alpha(\pp)$.
        \item When $\chi$ is trivial only: for any $\qq\mid\mid\nn$,  the Atkin-Lehner eigenvalue $\varepsilon(\qq)$.
    \end{enumerate}
\end{itemize}
Since we can only handle a finite amount of data in practice, we regard the input as an oracle which can supply the eigenvalue $\lambda_f(T)$ on demand for any principal $T$, and we will use this oracle a finite number of times.  The functions $\chi$ and $\varepsilon$ will be determined completely, having only finitely many values each, while we will only output the values $\alpha(\pp)$ for the finite set of primes whose norm is less than some bound and which do not divide the level.

By Theorem~\ref{thm:preprint-3.7}, the eigensystem $\lambda$ is only be determined up to unramified twist, but if we find one such eigensystem then we can find the others by twisting, so it suffices to find one.  In particular,  when the restriction of $\chi$ to $\CL[2]$ is trivial, so that $\chi$ is a square, we may choose $\chi$ to be trivial for the representative eigensystem in such a set of twists.   We proceed in stages, noting that the multiplicative relations satisfied by the eigenvalues mean that it suffices to determine $\chi$ fully, and $\alpha(\pp)$ for primes $\pp$ which do not divide the level. Recall that the quotient $\CL/\CL^2$ is called the \emph{genus group} of $K$; for an ideal $\aa$ we call its \emph{genus class} the image of its ideal class in the genus group.  We denote by $r_2$ the $2$-rank of the class group, so that $|\CL/\CL^2|=2^{r_2}$, and from classical genus theory we know that $r_2$ is one less than the number of primes dividing the discriminant of $K$.

\subsubsection{Outline of the procedure}
\begin{enumerate}
    \item Determine possible candidates for the character $\chi$, and fix one of these, taking $\chi$ to be trivial where possible.
    \item Determine the values $\alpha(\pp)$ for (a finite set of) prime ideals $\pp$ (not dividing $\nn$); each one will be processed using a method which depends on the ideal class and genus class of $\pp$:
    \begin{enumerate}
        \item $\alpha(\pp)$ for principal primes (those with trivial ideal class);
        \item $\alpha(\pp)$ for primes whose ideal class is square (those whose genus class is trivial);
        \item $\alpha(\pp)^2$ for primes in non-square classes (those whose genus class is non-trivial);
        \item lastly, use multiplicative relations to determine a consistent choice of signs for the $\alpha(\pp)$ with $\pp$ in each non-trivial genus class.  Different consistent choices of sign lead to eigensystems which are unramified quadratic twists of each other.
    \end{enumerate}
    \item If $\chi$ is trivial, determine the Atkin-Lehner eigenvalues $\varepsilon(\qq)$ for prime powers $\qq\mid\mid\nn$.
\end{enumerate}
The same procedure works in the case of odd class number, in which case steps 2(c,d) are not needed as the genus group is trivial.

\subsubsection{The procedure in detail}
\begin{enumerate}
    \item From  the known values of $\lambda_f(T_{\aa,\aa})$ for ideals $\aa$ coprime to $\nn$ with $\aa^2$ principal, we know the restriction of $\chi$ to the $2$-torsion subgroup $\CL[2]$.  We set $\chi$ to be an unramified character with this restriction, taking $\chi$ to be the trivial character when $\lambda_f(T_{\aa,\aa})=1$ for all such $\aa$.  This choice is always possible  in the odd class number case where $\CL[2]$ is trivial. 
    
    If our interest is restricted to eigensystems with trivial character, we only need to carry out this whole procedure on a homological eigenform $f$ for which $\lambda_f(T_{\aa,\aa})=1$ for all $\aa$ with $\aa^2$ principal, so we may first restrict to the subspace of $H_1(X_0(\mf{n});\QQ)$ on which  all principal operators $T_{\aa,\aa}$ act trivially.

    From now on we regard the character $\chi$ as fixed.
    \item We now consider the primes $\pp$ not dividing $\nn$ in order, determining $\alpha(\pp)$ for each:
    \begin{enumerate}
        \item If $\pp$ is principal, we have $\alpha(\pp)=\lambda_f(T_\pp)$ directly.
        \item If $\pp$ has square class, we find\footnote{This is always possible, since each ideal class contains ideals coprime to any fixed ideal, and even contains infinitely many primes.} an ideal $\aa$, coprime to $\pp\nn$, such that $\aa^2\pp$ is principal.   Using the principal eigenvalue $\lambda_f(T_{\aa,\aa}T_\pp) = \chi(\aa)\alpha(\pp)$,  we determine $\alpha(\pp) = \lambda_f(T_{\aa,\aa}T_\pp) / \chi(\aa)$.
        \item If $\pp$ does not have square class (equivalently, has non-trivial genus class), we use the relation
        \[
        \alpha(\pp)^2 = \alpha(\pp^2) + \n(\pp)\chi(\pp),
        \]
        together with the known principal eigenvalue $\lambda_f(T_{\aa,\aa}T_{\pp^2}) = \chi(\aa)\alpha(\pp^2)$,
        where $\aa$ is chosen in the inverse class to $\pp$ and coprime to $\nn$, to determine
        \[
        \alpha(\pp)^2 = \chi(\aa)^{-1}\lambda_f(T_{\aa,\aa}T_{\pp^2}) + \n(\pp)\chi(\pp) = \chi(\pp)(\lambda_f(T_{\aa,\aa}T_{\pp^2}) + \n(\pp)).
        \]
        This gives us the value of $\alpha(\pp)$ up to sign.  We need to choose the signs for all primes with non-trivial genus class in a systematic and consistent way. Of course, if $\alpha(\pp)^2=0$ then no choice is necessary.
        
        \item To deal with the primes in non-square classes, we maintain a table of pairs $(\aa,\alpha(\aa))$ for one ideal $\aa$ coprime to $\nn$ in each class $[\aa]$ in a subset of genus classes, representing the elements of $H/\CL^2$ for some subgroup $H$ with $\CL^2\subseteq H\subseteq \CL$, all with $\alpha(\aa)\not=0$.  Initially, $H=\CL^2$ and the table only contains $((1),1)$, where $(1)$ is the unit ideal. By construction, the ideals $\aa$ in the table will be square-free products of primes $\qq$ for which we have already decided (exactly, not just up to sign) the nonzero values $\alpha(\qq)$ from earlier choices.

        When we are processing a prime $\pp$ we first determine its ideal class and genus class.  If either of these is trivial, we use (2)(a) or (2)(b) to determine $\alpha(\pp)$.   Otherwise:
        \begin{itemize}
            \item find $\alpha(\pp)^2$ as in (2)(c); there is nothing more to do if this is zero.
            \item Otherwise, if the genus class is in our table, with entry $(\aa,\alpha(\aa))$, we set $\bb=\pp\aa$, which has trivial genus class and is square-free.  If $\bb$ is actually principal, we can compute $\alpha(\bb)=\lambda_f(T_\bb)$ and then $\alpha(\pp)=\alpha(\bb)/\alpha(\aa)$; otherwise we find an ideal $\dd$ with $\bb\dd^2$ principal, compute $\lambda_f(T_{\dd,\dd}T_\bb)$, and then $\alpha(\pp)=\lambda_f(T_{\dd,\dd}T_\bb) / (\alpha(\aa)\chi(\dd))$.
            \item Finally, suppose that the genus class of $\pp$ is not in the table.  Then we choose the sign of $\alpha(\pp)$ arbitrarily and update the table, doubling its size: for each entry $(\aa,\alpha(\aa))$ in the old table we add a new entry $(\aa\pp,\alpha(\aa)\alpha(\pp))$.
        \end{itemize}  
        Each time the table is extended, its size doubles; this happens at most $r_2$ times.  For eigensystems with no self-twist, this maximal size will be attained once we have found primes with nonzero eigenvalue $\alpha(\pp)$, whose genus classes generate the genus group.  The number of arbitrary choices made (after choosing $\chi$) is $r_2$, and the $2^{r_2}-1$ alternative eigensystems we would obtain by making different choices are precisely the nontrivial unramified quadratic twists of the original.
    \end{enumerate}
    \item We omit this step unless $\chi$ is trivial.
    
    To compute the eigenvalues $\varepsilon(\qq)$ of all Atkin-Lehner operators $W_\qq$ for the eigensystem,  the procedure is similar.  It suffices to do this for the prime power exact divisors $\qq\mid\mid\nn$.   Explicit $2\times2$ matrices for the principal operators $T$ required here are given in \cite[\S4.2]{cremona-preprint}, from which we obtain their eigenvalues $\lambda_f(T)$ as before.
    \begin{itemize}
        \item If $\qq$ is principal we compute $W_\qq$ directly to obtain $\varepsilon(\qq)=\lambda_f(W_\qq)$.
        \item Otherwise, if $\qq$ has square ideal class, we compute $T_{\aa,\aa}W_\qq$ for an ideal $\aa$ coprime to $\nn$ such that $\aa^2\qq$ is principal, and set $\varepsilon(\qq)=\lambda_f(T_{\aa,\aa}W_\qq)$.
        \item For $\qq$ with nontrivial genus class, we first observe that even if the eigensystem admits self-twist by a character $\psi$, we must have $\psi(\qq)=+1$, by Proposition \ref{prop:AL-restricts-self-twists}.  Hence there exists a prime $\pp$ in the inverse ideal class such that $\alpha(\pp)\not=0$, so we compute $T_\pp W_\qq$ and set $\varepsilon(\qq)=\lambda_f(T_\pp W_\qq)/\alpha(\pp)$.
\end{itemize}

\end{enumerate}

\begin{example}[Example \ref{ex:ur-twists-17} continued]
Consider again the field $K=\QQ(\sqrt{-17})$ with ring of integers $\OO=\ZZ[\omega]$ where $\omega=\sqrt{-17}$ and class group $\CL=\{1,c,c^2,c^3\}$ where $c=[(3,1+w)]$.  The character group is generated by $\chi_1$ of order $4$.

Given a homological eigensystem $\lambda_f$ at level $\nn$, we can determine the values of the associated character on $\CL[2]=\{1,c^2\}$ by evaluating $\lambda_f(T_{\aa,\aa})$ for one ideal $\aa$ coprime to $\nn$ with $[\aa]=c^2$.  If this value is $+1$ then the character of the full eigensystem is either trivial ($\chi_0$) or quadratic ($\chi_2$), and twisting by $\chi_1$ swaps these two choices, so we choose the trivial character and continue. On the other hand, if $\lambda_f(T_{\aa,\aa})=-1$ then the character is either $\chi_1$ or $\chi_3$, and we may choose either.

Having thus determined one of the two possibilities for the character $\chi$ of the full eigensystem, we determine the Hecke eigenvalues $\alpha(\pp)$, considering the primes not dividing $\nn$ one at a time.  For principal $\pp$ we have $\alpha(\pp)=\lambda_f(T_\pp)$, and if $[\pp]=c^2$ we have $\alpha(\pp)=\lambda_f(T_{\aa,\aa}T_\pp)/\chi(\aa)$, where $\aa$ is any ideal coprime to $\nn$ in class $c$ or $c^3$, so that $[\aa]^2=c^2$ and $[\aa^2\pp]=1$.

When we encounter a prime in ideal class $c$ or $c^3$ (that is, with non-trivial genus class), we compute $\alpha(\pp)^2$; if it is zero, we continue.  Suppose that we encounter one such prime with $\alpha(\pp)^2\not=0$. The first time this occurs, we pick one square root arbitrarily and set $\alpha(\pp)$ to this value, and we record the pair $(\pp,\alpha(\pp))$ in our table, which now has maximal size. From then on, given any other prime $\pp'$ in class $c$ or $c^3$, if $[\pp']=c^3$ then $\pp\pp'$ is principal and we can compute $\alpha(\pp')=\lambda_f(T_{\pp\pp'}) / \alpha(\pp)$, while if $[\pp']=c$ then $\alpha(\pp')=\lambda_f(T_{\pp,\pp}T_{\pp\pp'}) / (\chi(\pp)\alpha(\pp))$.

We may find that $\alpha(\pp)=0$ for all primes encountered in classes $c$ and $c^3$.  This is either because the eigensystem admits self-twist by $\chi_2$, or it may just be that we have not processed enough primes.  The former possibility can be ruled out if $\nn$ has a divisor $\qq\mid\mid\nn$ with $\chi_2(\qq)=-1$ (as in that case the eligible set $C(\nn)$ is empty), but otherwise we may have an ambiguous situation where the data so far considered from the homological eigensystem has not been enough to decide whether or not the eigenform admits self-twist.
\end{example}

\subsection{Galois actions}
There are two different Galois actions on the set of new Hecke eigenforms for the field $K$.   Let $\sigma$ generate $\Gal(K/\QQ)$; then it is easy to see that $\sigma$ induces an isomorphism from $S(\nn)^\ur$ to $S(\nn^\sigma)^\ur$, and also a map from $X_0(\nn)$ to $X_0(\nn^\sigma)$ which in turn induces an isomorphism on homology.  On eigensystems $\lambda=(\alpha,\chi)$ the action of sigma is simply $\lambda\mapsto\lambda^\sigma=(\alpha^\sigma,\chi^\sigma)$ where $\alpha^\sigma(\aa)=\alpha(\aa^\sigma)$ and similarly for $\chi$.  Hence, when the level $\nn$ is not Galois stable, we can easily find all eigensystems at the Galois  conjugate level $\nn^\sigma$ from those at level $\nn$ without having to repeat the whole algorithm; this saves time in practice when systematically finding all eigensystems at a range of levels.

Secondly, since $H_1(X_0(\nn),\QQ)$ is a vector space over $\QQ$, and all the principal Hecke operators are $\QQ$-linear, it follows that the eigenvalues for each operator are algebraic\footnote{In fact, the Hecke operators preserve integral homology, so the Hecke matrices will be integral if expressed with respect to an integral basis for the homology, so that the eigenvalues are algebraic integers; we note this fact, but do not use it in our computations.} and come in Galois conjugate sets.  Hence, in particular, the values of any homological eigensystem $\lambda_f$ are algebraic.  Since the homology is finite-dimensional, for each eigensystem the eigenvalues lie in a number field (of finite degree over $\QQ$) called the \emph{Hecke field} of the system, and for any $\tau\in\Gal(\overline{\QQ}/\QQ)$ the composition  $\tau\circ\lambda$ of an eigensystem with a Galois automorphism is another eigensystem at the same level.

    As first observed by Bygott, the field generated by the full eigensystem may be strictly larger than the field generated by the eigenvalues of principal operators.  There are examples in \cite{bygott} over $\QQ(\sqrt{-5})$ (with $\CLno=2$) where the principal eigenvalues are all rational, but the eigenvalues of nonprincipal operators generate a quadratic Hecke field.  One needs to take this possibility into account, for example, if seeking ``rational eigensystems'' with all eigenvalues rational.  Even if the character $\chi$ is trivial or quadratic (so only has  values $\pm1$) the process of computing $\alpha(\pp)$ for primes $\pp$ with nontrivial genus class involves extracting a square root, so we may find that a rational homological eigenform comes from a non-rational full Hecke eigenform.  The first author's code concentrates on rational eigensystems with trivial character (because of their connection to elliptic curves over $K$), and the implementation first finds rational homological eigenforms but discards some of these when some of the values $\alpha(\pp)^2$ are not rational squares.
   
\begin{example}[Example \ref{ex:ur-twists-17} continued]
    We keep the same notation as in Example \ref{ex:ur-twists-17}.  
    An example of this occurs over $\QQ(\sqrt{-17})$ at level $\nn=\pp=(2,1+\sqrt{-17})$.  The dimension of $H_1(X_0(\pp),\QQ)$ is $1$, so certainly all principal Hecke operators have rational eigenvalues, but the Hecke field is $\QQ(\sqrt{2})$.     In this example, there is only one homological eigensystem.   The character $\chi$ is trivial on $\CL[2]$ (as one checks that $T_{\aa,\aa}$ acts trivially for any ideal~$\aa$ of class $c^2$), hence the associated Hecke eigensystems have either the trivial character $\chi_0$ or the quadratic character $\chi_2$.  Each leads to two distinct Hecke eigensystems which are twists of each other, as well as being Galois conjugate.  Denoting one of those with trivial character by $\lambda$, the other is  $\lambda\otimes\chi_2 = \tau\circ\lambda$, where $\tau$ is the nontrivial automorphism of the Hecke field $\QQ(\sqrt{2})$. The Hecke eigensystems with character $\chi_2$ are the twists of these by the quartic character $\chi_1$, so have Hecke field $\QQ(\sqrt{2},\sqrt{-1})$.  See Example \ref{exm:hecke-sys-2a} for more details.
\end{example}

Examples such as this, where a Galois conjugate eigensystem is a twist of the original, are not typical.  Such an eigensystem is said to have \emph{Galois alignment} or \emph{inner twist}.  This occurs when there exists a non-trivial $\tau\in\Gal(\overline{\QQ}/\QQ)$ and character $\psi$ such that $\tau\circ\psi=\lambda^\sigma$; that is, if
\begin{equation}\label{eqn: inner-twist}
    \tau(\alpha(\pp)) = \psi(\pp) \alpha(\pp)
\end{equation}
for all prime ideals $\pp$ not dividing the level, where $\lambda=(\alpha,\chi)$.  When this occurs, we say that $\lambda$ \emph{admits inner twist} by $(\tau,\psi)$.  In the preceding example, the eigensystem with trivial character has inner twist by $(\tau,\chi_2)$ where (as above) $\tau$ generates $\Gal(\QQ(\sqrt{2})/\QQ)$.
\begin{lemma}
    If an eigensystem $\lambda=(\alpha,\chi)$ satisfies (\ref{eqn: inner-twist}) then it also satisfies
    \[
    \tau(\chi(\pp))=\psi(\pp)^2\chi(\pp),
    \]
    so that $\tau\circ\lambda = \lambda\otimes\psi$.
\end{lemma}
\begin{proof}
    This is the same proof as for classical modular forms: see \cite[Prop.~11.1.7(a)]{best2021computing}.  

    We need to show that if $\tau\circ\alpha=\psi\alpha$, then also $\tau\circ\chi=\psi^2\chi$.
    The character of $\tau\circ\lambda$ is $\tau\circ\chi$.  The multiplicative Hecke relations for $\lambda=(\alpha,\chi)$ and $\lambda\circ\psi=(\psi\alpha, \psi^2\chi)$ give
    \[
    \alpha(\pp)^2 - \alpha(\pp^2) = \chi(\pp)\n(\pp)
    \]
    and 
    \[
    \psi(\pp)^2\alpha(\pp)^2 - \psi(\pp^2)\alpha(\pp^2) = \psi(\pp)^2\chi(\pp)\n(\pp).
    \]
    Applying $\tau$ to the first equation and using the hypothesis gives
    \[
    \psi(\pp)^2\alpha(\pp)^2 - \psi(\pp)^2\alpha(\pp^2) = \tau(\chi(\pp))\n(\pp).
    \]
    Comparing the right-hand sides gives the result.
\end{proof}

\begin{example}[Example \ref{ex:ur-twists-17} continued]
We keep the same notation as in Example \ref{ex:ur-twists-17}.   There are two homological eigensystems which are Galois conjugate (over $\QQ(\sqrt{2})$ and twists by the quadratic character $\chi_2$.   Hence, there are four Hecke eigensystems, consisting of a single twist orbit, given in Table~\ref{tab:p2a}. Let $F_0$ denote the form in this orbit with trivial character, and set $F_j=F\otimes\chi_j=F\otimes\chi_1^j$ for $j=1,2,3$.  Then $F_0$ and $F_2$ have trivial character, have Hecke field $\QQ(\sqrt{2})$ and are quadratic twists of each other, while $F_1$ and $F_3$ have character $\chi_2$, Hecke field $\QQ(\sqrt{-2})$, and are also quadratic twists of each other.  Both $F_0$ and $F_2$ admit inner twist by $(\tau_1,\chi_2)$ where $\tau_1$ generates $\Gal(\QQ(\sqrt{2})/\QQ)$, while both $F_1$ and $F_2$ admit inner twist by $(\tau_2,\chi_2)$ where $\tau_2$ generates $\Gal(\QQ(\sqrt{-2})/\QQ)$.
\end{example}

\subsection{The Hecke field and principal subfield}
We end this section with some remarks on the relation between the Hecke field $k_F$ and its subfield, the principal Hecke field $k_f$. As already observed, $k_f$ is totally real, and $k_F$ contains $\QQ(\chi)$ where $\chi$ is the character of $F$.

First consider the case where $F$, and hence also $f$, has trivial character.  From our algorithmic procedure we see that the eigenvalues $\lambda_F(T)$ lie in the principal Hecke field $k_f$ for all operators $T$ in square classes, not just those which are principal.  (This is also true when $F$ has quadratic character, as it only uses the fact that the character values are rational.)   When the class number is odd, therefore, we have $k_F=k_f$. In general, for $T$ whose genus class is non-trivial, we only have $\lambda_F(T)^2\in k_f$, where the class of $\lambda_F(T)^2$ in $k_f^*/(k_f^*)^2$ only depends on the genus class of $T$; this is because if $T'$ has the same genus class as $T$ then $TT'$ has trivial genus class, so $\lambda(TT')\in k_f$.  Hence $k_F=k_f(\sqrt{r_1},\dots,\sqrt{r_m})$ for some $m\ge0$, where $r_1,\dots,r_m$ are multiplicatively independent modulo squares, and where $m$ is at most the $2$-rank of the class group; if $F$ admits a self-twist then $m$ is at most one less than this. For example, when the class group has cyclic $2$-primary part, $m\in\{0,1\}$; that is, either $k_F=k_f$ or $k_F=k_f(\sqrt{r})$ for some $r\in k_f^*$ which is not a square in $k_f$.  We will see examples of both these possibilities with $K=\QQ(\sqrt{-17})$ in the next section.

When the class number is odd, the character of $f$ is necessarily trivial, so that one lift, say $F_0$,  of the homological eigensystem also has trivial character with $k_{F_0}=k_f$, while all the other lifts are twists $F=F\otimes\chi$, so that $k_F=k_f(\chi)$, the extension of $k_f$ obtained by adjoining the values of the twisting character $\chi$.

The case where $f$ has non-trivial character is more complicated; rather than attempt a general description here, we will only give details in the case where the class group is cyclic of order $4$, as is relevant for the examples in the next section.

\section{Results and Tables for \texorpdfstring{$\QQ(\sqrt{-17})$}{Q(sqrt-17)}}
\label{sec:resultsandtables}

For this section, we fix the imaginary quadratic field $K=\QQ(\sqrt{-17})$ with a ring of integers $\OO = \ZZ[\omega]$, where $\omega = \sqrt{-17}$.  We use the \lmfdb\ labelling system (defined in \cite{CremonaPageSutherland}) to refer to integral ideals of $K$: each has a label of the form $N.i$ where $N$ is the ideal norm and $i$ its index in a well-defined ordering of the ideals of given norm. For example, if the rational prime $p$ splits in $K$, then the two primes above it are denoted $\pp_{p.1}$ and $\pp_{p.2}$.

The class group $\CL$ of $K$ is a cyclic group of order $4$. Let $c$ be a generator of $\CL$ such that $\pp_{3.1} \in c$, where $\pp_{3.1} = \langle 3, 4 + \omega \rangle$. The character group is also a cyclic group of order $4$.  Let $\chi_1$ be the character defined by $\chi_1(g) = i = \sqrt{-1}$,  and let $\chi_j = \chi_1^j$, so that $\chi_j(c^k) = i^{jk}$.  

\subsection{Summary of scope}
There are 305 ideals $\nn$ of norm less than or equal to 200. Of these, 244 have nontrivial homology $H_1(X_0(\nn),\QQ)$, and hence nontrivial spaces $S(\nn)^\ur$, and 104 have nontrivial newspaces $S(\nn)^{\ur,\new}$.  For each of these levels, we have computed all the Hecke eigensystems,  and the results are tabulated in Table~\ref{tab:newforms}. In this range of levels, we find 1100 new Hecke eigensystems, and hence $1100$ Bianchi newforms with unramified character, of which $8$ admit self-twist. There are 70 newforms that are potential base-change.  In this range,  $125$ eigensystems are rational ($62$ pairs of unramified twists, and one self-twist newform at level~$(8)=64.1$); in each of these cases we were able to find elliptic curves over $K$ of conductor equal to the level and with traces of Frobenius matching the Hecke eigenvalues at many small primes.

We also extended these computations to the range $200<\n(\nn) \le 1000$, but (to date) have limited our attention to newforms with trivial character and rational eigenvalues.  We find $651$ more (giving a total of $776$), which are tabulated in the \lmfdb\ at \url{https://www.lmfdb.org/ModularForm/GL2/ImaginaryQuadratic/2.0.68.1}.  Again, most of these newforms match elliptic curves over $K$, but there are four exceptions: a pair of quadratic twists at level $(6\omega-17)$ (of norm $901$), and another pair at the conjugate level.

\subsection{Tessellations}

To compute the homology $H_1(X_0(\mathfrak{n});\CC)$, we require a tessellation of $\hy_3$ with a $\go$ action. Here we summarize details about two such tessellations, as we use in our two independent implementations. 

The first tessellation comes from the theory of perfect Hermitian forms and their minimal vectors. For our field, there are $12$ perfect forms up to $\go$-equivalence.   The resulting tessellation of $\hy_3$ by ideal polyhedra has $12$ $\go$-orbits of polyhedra, consisting of $5$ tetrahedra, $2$ triangular prisms, $3$ square pyramids, one hexagonal cap, and one truncated tetrahedron. The $1$-homology is generated by the directed edges of these; in this case, we have $18$ $\go$-orbits of edges. The faces of the polyhedra give homology relations among these generators. In this case, there are $29$ $\GL{\OO}$-orbits of faces, including $22$ triangles, $6$ squares, and one hexagon, giving $3$-{}, $4$-, and $6$-term relations respectively.  Note that for a given level, we may not need all of the $\Gamma_0(\mathfrak{n})$-classes of these edges and faces to be included in our complex. There may be non-orientable edges and faces that come from the existence of orientation-reversing elements in the stabilizer groups of the polyhedra.

The second tessellation was constructed using a method based on Swan's algorithm. This has $13$ $\go$-orbits of polyhedra:  $6$ tetrahedra, $2$ cubes, $2$ triangular prisms, one hexagonal cap, one truncated tetrahedron, and one square pyramid. There are $14$ $\go$-orbits of edges which generate the $1$-homology, and $29$ orbits of faces giving edge-relations, similar to the Voronoi tessellation, in this case consisting of $20$ triangles, $8$ squares, and one hexagon.  

\subsection{Homological and full eigensystems and Hecke fields over \texorpdfstring{$\QQ(\sqrt{-17})$}{Q(sqrt(-17))}}

We discuss here in more detail for the field $\QQ(\sqrt{-17})$ the relation between the principal Hecke field $k_f$ and the full Hecke fields $k_F$ of the lifts of a homological newform $f$ to full eigenforms $F$.

At each level $\nn$, we decompose the rational new subspace $H$ of the homology $H_1(X_0(\nn), \QQ)$ as a module under the action of the (rational) Hecke algebra.  The irreducible submodules of $H$ correspond to Galois orbits of homological eigensystems.  Under the action of the principal involution operator $T_{\pp,\pp}$ (where $\pp$ is any prime in class $c^2$ coprime to the level), $H$ decomposes into two eigenspaces which we denote $H^+$ and $H^-$, with eigenvalue $+1$ and $-1$ respectively.  Each of these may further decompose under the action of the full principal Hecke algebra $\TT_1$,  and each irreducible component of dimension $d$ in $H^+$ or $H^-$ corresponds to a Galois orbit of $d$ homological newforms $f$, each with principal Hecke field $k_f$ a totally real number field of degree $d$; the Hecke fields of the Galois conjugates are isomorphic (conjugate) number fields.

Consider first the case of an irreducible component of $H^+$ of dimension $d$, with associated homological eigensystem $f$ of dimension $d$ and principal Hecke field $k_f$, of degree $d$ and totally real.  This lifts to full eigensystems $F_j$ for $j=0,1,2,3$:  we first lift to $F_0$ with trivial character $\chi_0$ and then twist, setting $F_j=F_0\otimes\chi_j$, so that $F_j$ has character $\chi_{2j}$.  That is, $F_0$ and its quadratic twist $F_2$ have trivial character $\chi_0$, while $F_1$ and its quadratic twist $F_3$ have character $\chi_2$.  In general, these four are distinct, but if $F_0$ has inner twist then (by definition) $F_2=F_0$ and also $F_1=F_3$.  

As for the Hecke fields, from our general discussion, we know that either $k_{F_0}=k_f$, or $k_{F_0}=k_f(\sqrt{r})$ for some non-square $r\in k_f$.  In the first case (which includes the case of self-twist),   the lifts of $f$ and its Galois conjugates to full eigensystems $F_0,F_2$ with characters $\chi_0$ fall into two Galois orbits (or just one in the self-twist case), each of size $d$.  In the second case, these lifts form a single Galois orbit of size $2d$, with full Hecke field $k_f(\sqrt{r})$ of degree $2d$.  Moreover, since the Hecke field of an eigensystem with trivial character is necessarily totally real, we see that the element $r\in k_f$ is totally positive.   The Hecke field of $F_2$ is $k_f(\sqrt{-1})$ in the first case, unless there is self-twist, when $k_{F_2}=k_f$, since the factor of $\sqrt{-1}$ only appears as a factor for eigenvalues of operators in classes $c$ and $c^3$, which are all zero in this case.  In the second case, we have $k_{F_2}=k_f(\sqrt{-r})$; there is one Galois orbit of size $2d$. (Since $r$ is totally positive, $-r$ is certainly not a square in the totally real field $k_f$.)  Note that in the second case all the full eigensystems have inner twist by $(\sigma, \chi_2)$ where $\sigma$ is the non-trivial automorphism of $\Gal(k_f(\sqrt{r})/k_f)$. Note also that all the eigenvalues are either elements of the principal field $k_f$, or products of such elements and either $\sqrt{r}$, $\sqrt{-1}$, or $\sqrt{-r}$.

Now consider the lifts of a $d$-dimensional homological eigensystem $f$ associated to an irreducible component of $H^-$.  Its lifts will have characters $\chi_1$ and $\chi_3$, two of each unless they have self-twist.   Since these characters are Galois conjugate (switching $i=\sqrt{-1}$ to $-i$), we only need consider the lifts with character~$\chi_1$, as those with character $\chi_3$ will be Galois conjugates of these, as well as being their twists.  Let $F$ be one lift to a full eigensystem with character $\chi_1$.  The eigenvalues for principal operators lie in $k_f$.  Those of operators in class $c^2$ are in $k_f(i)$, being multiples of elements of $k_f$ by $i$:  for example, if $\pp$ is a prime in class $c^2$, then taking $\aa$ in class $c$, we see that $T_{\pp}T_{\aa,\aa}$ is principal, so has eigenvalue in $k_f$, and the eigenvalue of $T_{\pp}$ is obtained by dividing by $\chi_1(\aa)=i$. So the full Hecke field certainly contains $k_f(i)$.  Next, for primes $\pp$ in class $c$ we see similarly that $\alpha(\pp)^2 = 2ri$ with $r\in k_f$ (the factor of $2$ is to simplify the formulas which follow), so $\alpha(\pp)=\pm\sqrt{2ri}=\pm\sqrt{r}(1+i)$.  So as before, we have two cases, according to whether $r$ is or is not a square in $k_f$. In the first case, the full Hecke field is $k_f(i)$, while in the second case it is the biquadratic extension $k_f(i, \sqrt{r})$. The lifts with character $\chi_1$ and $\chi_3$ form two Galois orbits of size $2d$ each in the first case, and a single orbit of size $4d$ in the second case;  in both cases, half of the Galois conjugates in each orbit have character $\chi_1$ and half $\chi_2$.  All the eigenvalues $\alpha(\pp)$ are products of an element of $k_f$ by either $1$ or $i$, or (in the second case) $\sqrt{r}$, or $\sqrt{r}(1+i)$.  Examples of both possibilities may be seen at level $\nn=(4)=16.1$: see Example~\ref{ex:level16.1} below.

\subsection{Newform Tables}\label{sec:new-tables}
In this section, we provide  summary tables of information about the newforms for levels of norm up to $200$.  

Table \ref{tab:newforms} shows the newspace dimensions for levels where the newspace is nontrivial, together with the dimensions of its decomposition into Hecke eigenspaces. Table \ref{tab:hecke-fields} gives the homological and full Hecke fields for each newform.

The columns of Table \ref{tab:newforms} are as follows:

\begin{itemize}\label{list:newform-table-details}
    \item $\mf{n}$,  $\overline{\mf{n}}$: \lmfdb\ labels of the level, and the conjugate level where this is different.
    \item $nd$: The dimension of the new Bianchi modular form space, $\dim S(\nn)^{\ur,\new}$; this is usually $4$  times $\dim H=\dim H^+ +\dim H^-$, but is less at level $64.1$ where there is a newform with self-twist; in the scope of this table, this only happens at level 64.1. See Example~\ref{ex:ur-twists-17}.
    \item $H^+$, $H^-$: The dimensions of irreducible components of $H^{\pm}$;   as a sequence of positive integers~$d$, each being the degree of the Hecke field $k_f$ of one Galois orbit of homological eigensystems $f$. The sum of these degrees (over both columns) is the dimension of $H$, the new part of the homology.    
    \item $\chi_0$: The sequence of dimensions of the Galois orbits of new Hecke eigensystems $F$ with trivial character $\chi_0$; again,  each integer is the degree of the full Hecke field $k_{F}$.  As explained above, for each entry $d$ in column $H^+$ we see either $d,d$ or $2d$ in this column, or just $d$ in the case of self-twist. The eigensystems with character $\chi_2$ are obtained from those of trivial character by quadratic twisting, with one degree $2d$ for each $d$ in column $H^+$.
    \item $\chi_1, \chi_3$: For each entry $d$ in the $H^-$ column, indicating a Galois orbit of $d$ homological eigensystems $f$ with principal Hecke field $k_f$ of degree $d$, there is either one entry $2d$ in this column (indicating one orbit of full eigensystems with self-twist and Hecke field $k_F=k_f(i)$), or two entries $2d$ (indicating two orbits of full eigensystems, each with Hecke field $k_F=k_f(i)$, quadratic twists of each other), or one entry $4d$ (indicating one orbit of full eigensystems $F$ with Hecke field $k_F=k_f(i,\sqrt{r})$ for some non-square $r\in k_f$).  In all cases, half the eigensystems in each Galois orbit have character $\chi_1$ and half $\chi_3$. 
\end{itemize}

See Examples~\ref{exm:hecke-sys-2a},   \ref{ex:level16.1}, \ref{ex:25.1}, and \ref{ex:modularity-example} for representative detailed examples. 

\begin{longtable}{L R C C C C} 
\caption{Newspace dimensions and decomposition into irreducible subspaces for the Hecke algebra.}
\label{tab:newforms}\\
\toprule
\mathfrak{n} , \mathfrak{\overline{n}} & nd & H^+ & H^- & \chi_0 & \chi_1, \chi_3 \\
\midrule
\endfirsthead
\caption{Newform table (continued)}\\
\toprule
\mathfrak{n} , \mathfrak{\overline{n}} & nd & H^+ & H^- & \chi_0 & \chi_1, \chi_3 \\
\midrule
\endhead
\bottomrule
\multicolumn{6}{c}{\emph{Continued on next page.}}
\endfoot
\bottomrule
\endlastfoot
2.1  & 4 & [ 1 ] & - & [ 2 ] & - \\
4.1  & 4 & [ 1 ] & - & [ 2 ] & - \\
7.1 , 7.2 & 4 & [ 1 ] & - & [ 1,1 ] & - \\
8.1  & 8 & [ 1,1 ] & - & [ 1,1,2 ] & - \\
9.1 , 9.3 & 4 & - & [ 1 ] & - & [ 2,2] \\
9.2  & 8 & [ 1,1 ] & - & [ 1,1,1,1 ] & - \\
12.1 , 12.2 & 4 & - & [ 1 ] & - & [ 4] \\
16.1  & 16 & - & [ 1,1,1,1 ] & - & [ 2,2,2,2,2,2,4] \\
17.1  & 4 & [ 1 ] & - & [ 1,1 ] & - \\
18.2  & 4 & [ 1 ] & - & [ 1,1 ] & - \\
21.1 , 21.4 & 4 & [ 1 ] & - & [ 1,1 ] & - \\
25.1  & 12 & [ 3 ] & - & [ 3,3 ] & - \\
26.1 , 26.2 & 4 & - & [ 1 ] & - & [ 2,2] \\
27.2 , 27.3 & 20 & [ 1 ] & [4] & [ 1,1 ] & [ 16] \\
34.1  & 4 & [ 1 ] & - & [ 1,1 ] & - \\
36.1 , 36.3 & 4 & [ 1 ] & - & [ 2 ] & - \\
36.2  & 4 & [ 1 ] & - & [ 1,1 ] & - \\
42.2 , 42.3 & 8 & [ 1 ] & [1] & [ 1,1 ] & [ 4] \\
48.1 , 48.2 & 4 & [ 1 ] & - & [ 1,1 ] & - \\
49.1 , 49.3 & 4 & [ 1 ] & - & [ 2 ] & - \\
49.2  & 20 & [ 5 ] & - & [ 5,5 ] & -  \\
50.1  & 12 & [ 1,1,1 ] & - & [ 1,1,1,1,1,1 ] & - \\
62.1 , 62.2 & 4 & [ 1 ] & - & [ 1,1 ] & - \\
63.1 , 63.6 & 8 & - & [ 1,1 ] & - & [ 4,4] \\
64.1  & 20 & [ 1,2 ] & [1,2] & [ 1,4 ] & [ 2,8] \\
66.1 , 66.4 & 4 & [ 1 ] & - & [ 1,1 ] & - \\
66.2 , 66.3 & 12 & - & [ 1,2 ] & - & [ 4,8] \\
68.1  & 8 & [ 2 ] & - & [ 2,2 ] & - \\
72.2  & 24 & [ 1,2,3 ] & - & [ 1,1,2,2,3,3 ] & - \\
78.2 , 78.3 & 20 & [ 1] & [4] & [ 1,1 ] & [ 16] \\
81.3  & 28 & [ 1,1,1 ] & [2,2] & [ 1,1,1,1,2 ] & [ 4,4,4,4] \\
93.2 , 93.3 & 8 & [ 1 ] & [1] & [ 1,1 ] & [ 4] \\
98.1 , 98.3 & 4 & [ 1 ] & - & [ 2 ] & - \\
98.2  & 28 & [ 1,1,2,3 ] & - & [ 1,1,1,1,2,2,3,3 ] & - \\
99.3 , 99.4 & 16 & - & [ 1,3 ] & - & [ 4,12] \\
100.1  & 16 & [ 1,3 ] & - & [ 1,1,3,3 ] & - \\
104.1 , 104.2 & 4 & [ 1 ] & - & [ 2 ] & - \\
106.1 , 106.2 & 4 & - & [ 1] & - & [2,2] \\
108.2 , 108.3 & 16 & [ 1,2 ] & [1] & [ 1,1,2,2 ] & [ 4] \\
112.1 , 112.2 & 4 & - & [ 1 ] & - & [ 4] \\
121.1 , 121.3 & 4 & - & [ 1 ] & - & [2,2] \\
121.2  & 36 & [ 1,8 ] & - &  [1,1,8,8 ] & - \\ 
126.1 , 126.6 & 12 & [ 1 ] & [1,1] & [ 1,1 ] & [ 4,4] \\
136.1  & 40 & [ 1,1,2,2 ] & [2,2] & [ 1,1,1,1,2,2,2,2 ] & [ 8,8] \\
138.2 , 138.3 & 4 & - & [ 1 ] & - & [ 4] \\
144.1 , 144.3 & 8 & [ 1 ] & [1] & [ 2 ] & [2,2] \\
144.2  & 40 & - & [ 1,1,1,1,1,2,3 ] & - & [ 4,4,4,4,4,8,12] \\
147.1 , 147.6 & 8 & [ 1,1 ] & - & [ 1,1,1,1 ] & - \\
150.1 , 150.2 & 4 & [ 1 ] & - & [ 1,1 ] & - \\
153.2  & 12 & [ 1,2 ] & - & [ 1,1,2,2 ] & - \\
159.2 , 159.3 & 4 & [ 1 ] & - & [ 1,1 ] & - \\
161.2 , 161.3 & 4 & [ 1 ] & - & [ 1,1 ] & - \\
162.1 , 162.5 & 16 & [ 2 ] & [ 1,1 ] & [ 4 ] & [ 4,4] \\
162.2 , 162.4 & 8 & [ 1 ] & [1] & [ 1,1 ] & [ 4] \\
162.3  & 16 & [ 1,1,1,1 ] & - & [ 1,1,2,2,2 ] & - \\
168.1 , 168.4 & 8 & [ 1 ] & [1] & [ 1,1 ] & [ 4] \\
169.2  & 36 & [ 1,2,3,3 ] & - & [ 1,1,4,6,6 ] & - \\
178.1 , 178.2 & 16 & [1] & [ 1,1,1 ] & [ 1,1 ] & [ 2,2,2,2,2,2] \\
186.1 , 186.4 & 4 & - & [ 1 ] & - & [ 4] \\
189.4 , 189.5 & 8 & - & [ 1,1 ] & - & [ 4,4] \\
196.1 , 196.3 & 4 & - & [ 1 ] & - & [2,2] \\
196.2  & 16 & [4] & - & [ 4,4 ] & - \\
198.2 , 198.5 & 20 & [3] & [1,1] & [ 3,3 ] & [4,4] \\
198.3 , 198.4 & 8 & - & [ 1,1 ] & - & [4,4] \\
200.1 & 72 & [ 1,1,4,4 ] & [1,1,3   ,3] & [ 1,1,1,1,4,4,4,4 ] & [2,2,2,2,12,12] \\
\end{longtable}

In Table \ref{tab:hecke-fields} we give the Hecke fields for the Bianchi newforms with trivial character.  We only include one of each pair of conjugate levels $\nn$, $\overline{\nn}$ as the entries would be the same for each, and we also only include one of each pair of unramified twists.   Where there is more than one newform (up to quadratic twist) at the same level, we number these $1,2,3,\dots$.  For homological Hecke fields $k_f$ of degree greater than $2$ we give the \lmfdb\ label of the field where possible; the degree $8$ field $\QQ(a)$ which occurs at level $121.2=(11)$ is not in the \lmfdb,  and has defining polynomial $f_8(x) = x^{8} -  x^{7} - 10 x^{6} + 8 x^{5} + 28 x^{4} - 18 x^{3} - 20 x^{2} + 6 x +  4$.  For newforms where $k_F\not=k_f$, we give $k_F$ in the form $k_f(\sqrt{r})$, where $r$ is a non-square in $k_f$.

\begin{longtable}{R R | C C || R R | C C}
\caption{Hecke fields for newforms with trivial character}
\label{tab:hecke-fields}\\
\toprule
\mathfrak{n} & \# & k_f & k_F & \mathfrak{n} & \# & k_f & k_F  \\
\midrule
\endfirsthead
\caption{Hecke fields (continued)}\\
\toprule
\mathfrak{n} & \# & k_f & k_F & \mathfrak{n} & \# & k_f & k_F  \\
\midrule
\endhead
\bottomrule
\multicolumn{8}{c}{\emph{Continued on next page.}}
\endfoot
\bottomrule
\endlastfoot
2.1  & 1 & \QQ & \Qroot{2} &
100.1 & 1 & \QQ & \QQ \\
4.1  & 1 & \QQ & \Qroot{2} &
100.1 & 2 & \lmfdbfield{3.3.1524.1} & \lmfdbfield{3.3.1524.1} \\
7.1  & 1 & \QQ & \QQ &
104.1 & 1 & \QQ & \Qroot{2} \\
8.1  & 1 & \QQ & \QQ &
108.2 & 1 & \QQ & \QQ \\
8.1  & 2 & \QQ & \Qroot{2} &
108.2 & 2 & \Qroot{15} & \Qroot{15} \\
9.2  & 1 & \QQ & \QQ &
121.2 & 1 & \QQ & \QQ \\
9.2  & 2 & \QQ & \QQ &
121.2 & 2 & \QQ(a), f_8(a)=0 & \QQ(a) \\
17.1 & 1 & \QQ & \QQ &
126.1 & 1 & \QQ & \QQ \\
18.2 & 1 & \QQ & \QQ &
136.1 & 1 & \QQ & \QQ \\
21.1 & 1 & \QQ & \QQ &
136.1 & 2 & \QQ & \QQ \\
25.1 & 1 & \lmfdbfield{3.3.148.1} & \lmfdbfield{3.3.148.1} &
136.1 & 3 & \Qroot{5} & \Qroot{5} \\
27.2 & 1 & \QQ & \QQ &
136.1 & 4 & \Qroot{3} & \Qroot{3} \\
34.1 & 1 & \QQ & \QQ &
144.1 & 1 & \QQ & \Qroot{2} \\
36.1 & 1 & \QQ & \Qroot{2} &
147.1 & 1 & \QQ & \QQ \\
36.2 & 1 & \QQ & \QQ &
147.1 & 2 & \QQ & \QQ \\
42.2 & 1 & \QQ & \QQ &
150.1 & 1 & \QQ & \QQ \\
48.1 & 1 & \QQ & \QQ &
153.2 & 1 & \QQ & \QQ \\
49.1 & 1 & \QQ & \Qroot{2} &
153.2 & 2 & \Qroot{17} & \Qroot{17} \\
49.2 & 1 & \lmfdbfield{5.5.240133.1} & \lmfdbfield{5.5.240133.1} &
159.2 & 1 & \QQ & \QQ \\
50.1 & 1 & \QQ & \QQ &
161.2 & 1 & \QQ & \QQ \\
50.1 & 2 & \QQ & \QQ &
162.1 & 1 & \Qroot{3} & \QQ(\sqrt{2},\sqrt{3}) \\
50.1 & 3 & \QQ & \QQ &
162.2 & 1 & \QQ & \QQ \\
62.2 & 1 & \QQ & \QQ &
162.3 & 1 & \QQ & \QQ \\
64.1 & 1 & \QQ & \QQ &
162.3 & 2 & \QQ & \Qroot{2} \\
64.1 & 2 & \Qroot{2} & \QQ(\sqrt{2+\sqrt{2}}) &
162.3 & 3 & \QQ & \Qroot{2} \\
66.1 & 1 & \QQ & \QQ &
162.3 & 4 & \QQ & \Qroot{2} \\
68.1 & 1 & \Qroot{3} & \Qroot{3} &
168.1 & 1 & \QQ & \QQ \\
72.2 & 1 & \QQ & \QQ &
169.2 & 1 & \QQ & \QQ \\
72.2 & 2 & \Qroot{33} & \Qroot{33} &
169.2 & 2 & \Qroot{17} & \QQ(\sqrt{(5+\sqrt{17})/2}) \\
72.2 & 3 & \lmfdbfield{3.3.316.1} & \lmfdbfield{3.3.316.1} &
169.2 & 3 & \QQ(a)=\lmfdbfield{3.3.621.1} & \QQ(\sqrt{a^2-2a-1}) \\
78.2 & 1 & \QQ & \QQ &
169.2 & 4 & \QQ(a)=\lmfdbfield{3.3.229.1} & \QQ(\sqrt{5-a^2}) \\
81.3 & 1 & \QQ & \QQ &
178.1 & 1 & \QQ & \QQ \\
81.3 & 2 & \QQ & \QQ &
196.2 & 1 & \lmfdbfield{4.4.1957.1} & \lmfdbfield{4.4.1957.1} \\
81.3 & 3 & \QQ & \Qroot{17} &
198.2 & 1 & \lmfdbfield{3.3.788.1} & \lmfdbfield{3.3.788.1} \\
93.2 & 1 & \QQ & \QQ &
200.1 & 1 & \QQ & \QQ \\
98.1 & 1 & \QQ & \Qroot{2} &
200.1 & 2 & \QQ & \QQ \\
98.2 & 1 & \QQ & \QQ &
200.1 & 3 & \lmfdbfield{4.4.23252.1} & \lmfdbfield{4.4.23252.1} \\
98.2 & 2 & \QQ & \QQ &
200.1 & 4 & \lmfdbfield{4.4.92692.1} & \lmfdbfield{4.4.92692.1} \\
98.2 & 3 & \Qroot{3} & \Qroot{3} &&&&\\
98.2 & 4 & \lmfdbfield{3.3.148.1} & \lmfdbfield{3.3.148.1} &&&&\\
\end{longtable}

\pagebreak
\subsection{Detailed examples}


\begin{example}\label{exm:hecke-sys-2a}
The first non-trivial Hecke eigensystem occurs at level $\pp_{2.1}=\langle 2,1 + \omega \rangle$. The excerpt from Table~\ref{tab:newforms} corresponding to this level is given below. 
\[
  \begin{array}{c c c c c c }
  \toprule
\mathfrak{n} & nd & H^+ & H^- & \chi_0 & \chi_1, \chi_3 \\ \hline
2.1 & 4 & [ 1 ] & - & [ 2 ] & - \\
  \bottomrule
  \end{array}
  \]
At this level, there is only one homological eigenform $f$. Since the homology is one-dimensional, the principal eigenvalues $\lambda_f(T)$ of $f$ are rational, so the principal Hecke field $k_f$ is $\QQ$.  

Let $\pp = \mf{p}_{13.1} = \langle 13,3 + \omega \rangle$.  Then $\mf{p} \in \CL^2$, and the eigenvalue of the principal Hecke operator $T_{\pp,\pp}$ is $1$. Hence $H=H^+$, while $H^-=0$.   Therefore, the character orbit of $f$ is $\set{\chi_0,\chi_2}$. We fix the character to be the trivial character $\chi=\chi_0$ to identify the full Hecke eigensystem $(\alpha,\chi)$ for one newform $F_0$.

The full eigensystems $F_0,F_2$ in the twist orbit which have trivial character have full Hecke field $\QQ(\sqrt{2})$ while $F_1,F_3$ have character $\chi_2$ and full Hecke field $\QQ(\sqrt{-2})$, as indicated by the entry $[2]$ in the $\chi_0$ column. 

The first primes not dividing the level are the two primes above 3, $\pp_{3.1}=\langle 3,4 + \omega \rangle \in c$ and $\pp_{3.2} = \langle 3, 2+\omega  \rangle\in c^3$. The operator $T = T_{\pp_{3.1},\pp_{3.1}}T_{\pp_{3.1}^2}$ is principal and has eigenvalue $\lambda_f(T)=5$, so that $\alpha(\pp_{3.1}^2)=5$ and hence $\alpha(\pp_{3.1})^2 = 5+\chi(\pp_{3.1})\n(\pp_{3.1}) = 8$.  As this is nonzero, we may arbitrarily fix a sign for the square root, taking $\alpha(\pp_{3.1})=2\sqrt{2}$.  At this point, we have made two choices and thus determined the full eigensystem $(\alpha,\chi)$ uniquely; there is no self-twist, so the twist orbit of $F_0$ consists of four distinct newforms $F_i=F\otimes\chi_i$ for $0\le i\le 3$.

To compute $\alpha(\pp_{3.2})$ we compute the principal operator $T = T_{\pp_{3.1}\pp_{3.2}}$, which has eigenvalue $\lambda_f(T) = -8$, obtaining $\alpha(\pp_{3.2})=-8/(2\sqrt{2}) = -2\sqrt{2}$.    Similarly, for all other primes $\pp$ in class $c^3$ we can obtain $\alpha(\pp)$ from the eigenvalue of the principal operator $T = T_{\pp\pp_{3.1}}$ as $\alpha(\pp) = \lambda_f(T)/\alpha(\pp_{3.1}) = \lambda_f(T)/(2\sqrt{2})$.  For primes $\pp$ in class $c$ we use the principal operator $T = T_{\pp\pp_{3.2}}$ to obtain $\alpha(\pp) = \lambda_f(T)/\alpha(\pp_{3.2}) = \lambda_f(T)/(-2\sqrt{2})$.   Since all the principal operators have rational eigenvalues, we see that all the non-principal ones have eigenvalues in $\QQ(\sqrt{2})$; in fact they are all rational multiples of $\sqrt{2}$.

For primes in class $c^2$, we use the principal operator $T=T_{\pp,\pp}T_{\pp}$ and compute $\alpha(\pp)=\lambda_f(T)$.  For example,  for $\pp=\pp_{13.1}$ we find that $\lambda_f(T_{\pp,\pp}T_{\pp}) = -2$, giving $\alpha(\pp)=-2$.

Since the level is prime, there is just one Atkin-Lehner eigenvalue $\varepsilon_\qq$ for $\qq=\nn$.  Taking $\aa=\pp_{3.1}$ so that $\aa^2\qq$ is principal, we find that $\lambda_f(T_{\aa,\aa}W_\qq) = -1$.  Hence $\varepsilon_\qq(F_0) = -1$ (and $\varepsilon_{\qq}(F_2)=-1$ also).

In summary, we have computed the eigenvalues of one full eigensystem $(\alpha,\chi)$ associated to one newform $F_0$ using principal Hecke operators.  The full Hecke field is $k_{F_0}=\QQ(\sqrt{2})$.   As this eigensystem does not admit self-twist, its twist orbit has size four;  twisting this by all characters of the class group gives the four Hecke eigensystems given in Table \ref{tab:p2a}, where $\alpha = \sqrt{2}$ and $\beta = \sqrt{-2}$. The newforms $F_0$ and $F_2$ have trivial character and Hecke field $\QQ(\sqrt{2})$, while $F_1$ and $F_2$ have character $\chi_2$ and Hecke field $\QQ(\sqrt{-2})$.
 
\begin{table}[ht]
\caption{Hecke eigensystems at level $\pp_{2.1}$. ($\alpha = \sqrt{2}, \beta = \sqrt{-2} = i\sqrt{2}$)} 
\label{tab:p2a}
\begin{tabular}{|C|RRRRRRRRRRRR|}
\hline
\pp
  & \pp_{3.1} & \pp_{3.2} & \pp_{7.1} & \pp_{7.2} & \pp_{11.1} & \pp_{11.3} & \pp_{13.1} & \pp_{13.2} & \pp_{17.1} & \pp_{23.1} & \pp_{23.2} & \pp_{25.1} \\ \hline
[\pp] & c^3 & c & c^3 & c & c & c^3 & c^2 & c^2 & 1 & c & c^3 & 1 \\ \hline
F_0  & 2\alpha & -2\alpha & 0 & 0 & 2\alpha & -2\alpha & -2 & -2 & -6 & -4\alpha & 4\alpha & 2 \\
F_1=F\otimes \chi_1 &  2\beta & 2\beta & 0 & 0 & -2\beta & -2\beta & 2 & 2 & -6 & 4\beta & 4\beta & 2 \\
F_2=F\otimes \chi_2 &  2\alpha & -2\alpha & 0 & 0 & 2\alpha & -2\alpha & -2 & -2 & -6 & 4\alpha & -4\alpha & 2 \\
F_3=F\otimes \chi_3 &  -2\beta & -2\beta & 0 & 0 & 2\beta & 2\beta & 2 & 2 & -6 & -4\beta & -4\beta & 2 \\ \hline
\end{tabular}
\end{table}

Note that the newforms $F_1$ and $F_3$ have the same eigenvalues on conjugate primes, suggesting that they are potentially base change of a classical newform.   In fact this must be the case, since the level $\nn=\pp_{2.1}$ is stable under Galois conjugation, so for every newform $G\in S(\nn)^{\ur}$ the newform $G^{\sigma}$ is also in $S(\nn)^\ur$ (where, as above, $\sigma$ denotes the generator of $\Gal(K/\QQ)$).  Inspection of Table \ref{tab:p2a} suffices to see how $\sigma$ acts on the four newforms:  we see that $F_0^\sigma=F_2$ (and vice versa), while both $F_1$ and $F_3$ are fixed by~$\sigma$.   It follows, by a result of G\'erardin and Labesse \cite{GerardinLabesse} that the latter two newforms are indeed base-change. The eigensystem $F_1$ agrees exactly with a classical newform $g$ at level $34$ with a quadratic Dirichlet character $\psi$ of conductor $34$ \href{https://www.lmfdb.org/ModularForm/GL2/Q/holomorphic/34/2/b/a/}{(34.2.b.a)}, on split primes. On inert primes, we have the following relationship. 
\[
\begin{array}{|c|cccc|cc}\hline
     p &        5 & 19 & 43&59  \\ \hline
     \psi(p) & -1 & 1 & 1 & 1\\
     a_p     &-2\sqrt{-2}& -4 & -4&12 \\ \hline
     a_{\pp} = a_p^2-2\psi(p)p& 2 & -22 & -70 & 26  \\ \hline 
\end{array}
\]

Moreover, we observed that in the scope of our computation, the L-function of $F_0$ matches the L-function of the two Hilbert cusp form $h_1,h_2$ with \lmfdb  labels
\href{https://www.lmfdb.org/ModularForm/GL2/TotallyReal/2.2.17.1/holomorphic/2.2.17.1-32.3-a}{2.2.17.1-32.3-a} and \href{https://www.lmfdb.org/ModularForm/GL2/TotallyReal/2.2.17.1/holomorphic/2.2.17.1-32.4-a}{2.2.17.1-32.4-a}. Furthermore, $h_1$ and $h_2$ match the isogeny classes of two elliptic curves of conductor $32$ over $\QQ(\sqrt{17})$ \href{https://www.lmfdb.org/EllipticCurve/2.2.17.1/32.3/a/}{(32.3a and }  
\href{https://www.lmfdb.org/EllipticCurve/2.2.17.1/32.4/a/}{32.4a)}. 

On the other hand, the base change of $g$ to $\QQ(\sqrt{17})$ has rational coefficients because it has inner twist by the Dirichlet character \href{https://www.lmfdb.org/Character/Dirichlet/17/b}{17.b}. By modularity results for Hilbert modular forms, we can recover the same two isogeny classes of elliptic curves as above. 

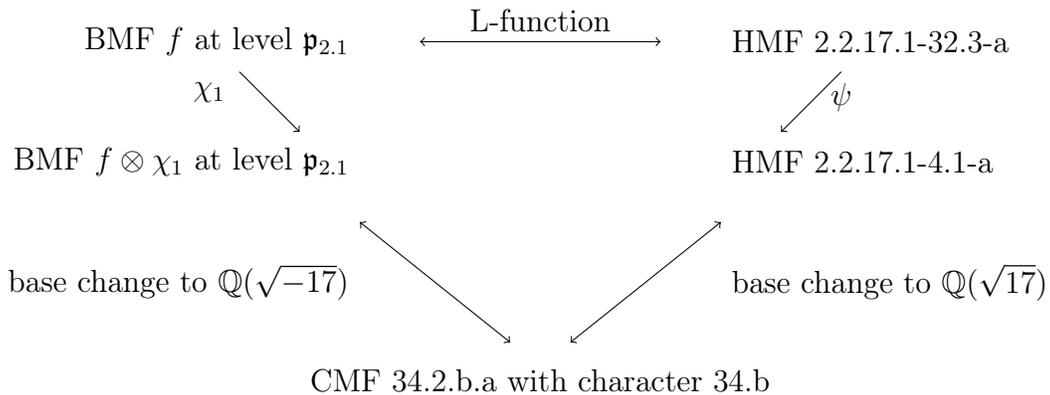
\begin{figure}[ht!]
    \centering
    \begin{tikzpicture}[scale=0.8]
        \draw (0,0) node [above] {CMF \href{https://www.lmfdb.org/ModularForm/GL2/Q/holomorphic/34/2/b/a/}{34.2.b.a} with character \href{https://www.lmfdb.org/Character/Dirichlet/34/b}{34.b}};

        \draw (-3,2)node[left]{base change to $\QQ(\sqrt{-17})$};
        \draw (-3,4) node [left] {BMF $f\otimes \chi_1$ at level $\pp_{2.1}$} ;
        \draw[<-](-4,4.5) -- (-5,5.5);
        \draw (-5.5,5.5) node [below] {$\chi_1$};

        \draw (3,2)node[right]{base change to $\QQ(\sqrt{17})$};
        \draw (3,4)node[right]{HMF \href{https://www.lmfdb.org/ModularForm/GL2/TotallyReal/2.2.17.1/holomorphic/2.2.17.1-4.1-a}{2.2.17.1-4.1-a}};
        \draw[<-](4,4.5) -- (5,5.5) node [below] {$\psi$};

        \draw[<->] (-0.5,1) -- (-3,3);
        \draw[<->] (0.5,1) -- (3,3);

        \draw (-3,6)node[left]{BMF $f$ at level $\pp_{2.1}$};

        \draw (3,6)node[right]{HMF \href{https://www.lmfdb.org/ModularForm/GL2/TotallyReal/2.2.17.1/holomorphic/2.2.17.1-32.3-a}{2.2.17.1-32.3-a}};
        \draw [<->] (-2,6) -- (2,6);
        \draw (0,6) node [above] {L-function};
    \end{tikzpicture}
    \caption{Observed connections at level $\pp_{2.1}$}
\end{figure}

\end{example}


\begin{example}\label{ex:level16.1} 
The row in the newform tables corresponding to the level $16.1$ from Table~\ref{list:newform-table-details} is given below.
\[
  \begin{array}{cccccc}
  \toprule
\mathfrak{n} &  nd & H^+ & H^- & \chi_0 & \chi_1, \chi_3 \\ \hline
16.1 & 16 & - & [ 1,1,1,1 ] & - & [ 2,2,2,2,2,2,4] \\
\bottomrule
  \end{array}
\]
The newspace $H$ has dimension $4$, and we find $4$ new homological eigensystems $f$, each with rational coefficient field, $k_f=\QQ$. This level does not contain any self-twist systems; therefore, the total newform dimension is $4 \cdot 4 = 16$. The first prime not dividing the level in class $c^2$ is $\pp=\pp_{13.1}$, and the principal operator $T_{\pp,\pp}$ acts as $-1$ on the whole newspace~$H$.  To see that $H$ splits into four one-dimensional irreducible components under the action of $\TT_1$, it suffices to consider the principal operator $T_{\pp}T_{\pp,\pp}$ with $\pp=\pp_{3.1}$ (in class $c$), which has distinct eigenvalues $-5,-2,2,3$.  This explains the entries of $0$ and $[1,1,1,1]$ in the $H^+$ and $H^-$ columns respectively.

Each homological eigensystem lifts to two full eigensystems with character $\chi_1$ and two with the conjugate character $\chi_3$. One of them lifts to a single Galois orbit of size $4$ and Hecke field $k_F=\QQ(i,\sqrt{2}) = \QQ(\zeta)$ where $\zeta^4=-1$, consisting of two newforms with each character.  The other three homological eigensystems lift to full eigensystems with Hecke field $k_F=\QQ(i)$.

The eigenvalues of one Hecke eigensystem in each Galois orbit, for the first few primes are given in Table \ref{tab:16.1}.  Newform $F_1$ has the degree $4$ Hecke field $\QQ(i,\sqrt{2})$, and is conjugate to its quadratic twist.  The other six Galois orbits are in pairs of quadratic twists $\{F_2,F_3\}$, $\{F_4,F_5\}$, $\{F_6,F_7\}$; the eigenvalues in the table are for the conjugate with character $\chi_1$.

\begin{table}[ht]  
\caption{Hecke eigensystems at level $16.1 = \pp_{2.1}^4$.}  
\label{tab:16.1}
\resizebox{\textwidth}{!}{
\begin{tabular}{|C|CCCCCCCCC|}
\hline
 \pp & \pp_{3.1} & \pp_{3.2} & \pp_{7.1} & \pp_{7.2} & \pp_{11.1} & \pp_{11.2} & \pp_{13.1} & \pp_{13.2} & \pp_{17.1} \\ 
\hline
 F_1 & \sqrt{2}(1+i) & -\sqrt{2}(1-i) & -\sqrt{2}(1+i) & \sqrt{2}(1-i) & \sqrt{2}(1-i) & -\sqrt{2}(1+i) & -2i & 2i & -2  \\
 F_2 & 1+i & -1+i & 3+3i & -3+3i & -1+i & 1+i & -4i & 4i & -6 \\
 F_3 & -1-i & 1-i & -3-3i & 3-3i & 1-i & -1-i & -4i & 4i & -6 \\
 F_4 & 2+2i & -2+2i & 0 & 0 & -2+2i & 2+2i & 2i & -2i & 6  \\
 F_5 & -2-2i & 2-2i & 0 & 0 & 2-2i & -2-2i & -2i & 2i & 6  \\ 
 F_6 & 0 & 0 & 2+2i & -2+2i & 4-4i & -4-4i & 2i & -2i & 6  \\
 F_7 & 0 & 0 & -2-2i & 2-2i & -4+4i & 4+4i & -2i & 2i & 6  \\
\hline
\end{tabular}
}
\end{table}

\end{example}


\begin{example} \label{ex:25.1}
We describe the eigensystems at level $\nn = 25.1 = 5\OO$,  the prime ideal above $5$.
 Let $L=\QQ(a)$ be the degree $3$ field defined (up to conjugation) by the polynomial $f(x)= x^3 - x^2 - 3x + 1$, with LMFBD label \href{https://www.lmfdb.org/NumberField/3.3.148.1}{3.3.148.1}.
There is a homological eigensystem $f$ at level $\nn$ with trivial character and principal Hecke field $k_f=L$, in a Galois orbit of three such eigensystems.  The full Hecke field is also $L$.  The first few eigenvalues $a_{\pp}$ of one full eigensystem with trivial character are given in Table~\ref{tab:25.1}.  
\begin{table}[ht]  
\caption{One Hecke eigensystem at level $25.1 = \pp_{5.1}^2$. ($a$ denotes a root of $f(x)$.)}  
\label{tab:25.1}
\begin{tabular}{|C|CCCCCCC|}
\hline
\pp
  & \pp_{2.1} & \pp_{3.1} & \pp_{3.2} & \pp_{7.1} & \pp_{7.2} & \pp_{11.1} & \pp_{11.2} \\ 
  \hline 
a_{\pp}
  & a         & a^2-a-2   & -a^2+a+2  & a^2+a-2   & -a^2-a+2  & -a-1       & a+1\\
\hline
\end{tabular}
\vskip10pt
\begin{tabular}{|C|CCCCC|}
\hline
\pp 
  & \pp_{13.1} & \pp_{13.2} & \pp_{17.1} & \pp_{23.1} & \pp_{23.2} \\ \hline 
a_{\pp}
& a^2-2a-3 & a^2-2a-3 & -4a^2+2a+8 & a^2+3a-6 & -a^2-3a+6  \\
\hline
\end{tabular}
\end{table}
Twisting this orbit by the quadratic character $\chi_2$ gives another Galois orbit with trivial character and Hecke field $L$ at level $\nn$.  Twisting  by  $\chi_1$ or $\chi_3$ gives eigensystems with character $\chi_2$ and Hecke field $L(i)$ at level $\nn$, in a single Galois orbit of size $6$.   All the eigensystems with trivial character have Atkin-Lehner eigenvalue $+1$.

The row in Table \ref{tab:newforms} corresponding to the level $25.1$ is given below, where we have added a column to show the size of the single Galois orbit of full eigensystems with character $\chi_2$:
\[
  \begin{array}{ccccccc}
    \toprule
  \mathfrak{n} & nd & H^+ & H^- & \chi_0 & \chi_2 & \chi_1, \chi_3 \\
  \midrule
  25.1 & 12 & [ 3 ] & - & [ 3,3 ] & [ 6 ] & - \\ 
\bottomrule
  \end{array}
  \]

\end{example}


\begin{example}
\label{ex:modularity-example}
The first rational Hecke eigensystems occur at each of the two conjugate levels $\pp_{7.1} = \langle7,w+2\rangle$ and $\pp_{7.2} = \langle 7, 5 + \omega \rangle$. For this example we have chosen a newform $F$ at level $\pp_{7.2}$, which may be found in the \lmfdb\ where its label is \lmfdbbmf{2.0.68.1}{7.2}{a}.   We show how to use the algorithm given in \cite{modularityexamplesDGP}, based on the Serre-Faltings-Livn\'e method, to prove the modularity of an associated elliptic curve $E$ defined over $K$ with conductor $\pp_{7.2}$, which has \lmfdb\ label \lmfdbecnf{2.0.68.1}{7.2}{a}{2} and equation
\[
E: y^2 + w xy = x^3 +(w+1) x^2 -(w+8) x.
\]
(Note that this curve $E$ does not have a global minimal model; the model we use here is non-minimal at the prime $\pp_{3.1}=\langle3,w+2\rangle$.)

At level $\pp_{7,2}$, the homological dimension is $1$, so there is one homological eigenform, which gives rise to two full eigensystems with trivial character and labels \lmfdbbmf{2.0.68.1}{7.2}{a} and \lmfdbbmf{2.0.68.1}{7.2}{b}, these being quadratic twists of each other. Table \ref{tab:7.2} gives the first few eigenvalues of $T_{\pp}$ (for primes $\pp$ not dividing the level) for the first of these.
\begin{table}[ht]
\caption{The Hecke eigensystem at level $\pp_{7.2}$.}
\label{tab:7.2}
\begin{tabular}{|C|CCCCCCCCCCCC|}
\hline
 \pp & \pp_{2.1} & \pp_{3.1} & \pp_{3.2} & \pp_{7.1} & \pp_{11.1} & \pp_{11.2} & \pp_{13.1} & \pp_{13.2} & \pp_{17.1} & \pp_{23.1} & \pp_{23.2} & \pp_{25.1} \\ \hline
a_{\pp}(F) & -1 & -2 & -2 & 0 & -6 & 2 & -2 & 6 & -2 & 0 & 8 & -2  \\
  \hline 
  \end{tabular}
\end{table}
One readily checks that the elliptic curve $E$ has the same values for the trace of Frobenius $a_{\pp}(E)$ at each of these primes, and one can in principle carry out the same check for any finite number of primes.  In order to show that $E$ is modular, it is necessary in principal to establish that $E$ and the newform $F$ have equal $a_{\pp}$ for all (or at least all but finitely many) primes $\pp$, as this implies that their $\ell$-adic Galois representations are isomorphic for one, and hence all primes $\ell$, and that their $L$-functions are equal.

We now explain how to carry out the algorithm of \cite[Algorithm~2.2]{modularityexamplesDGP}.   This involves determining a finite set of primes $\pp$ such that the equality of the trace of Frobenius $a_{\pp}(E)$ and the Hecke eigenvalue $a_{\pp}(F)$ for the primes in this set is sufficient to imply their equality for all primes.  We refer the reader to \cite[Algorithm~2.2]{modularityexamplesDGP} for more details of the algorithm.

\begin{enumerate}
    \item First we chose two prime ideals $\pp = \pp_{2.1}$ and $\pp_{3.1}$.  These have nonzero trace $a_{\pp}$, and in each case the prime $2$ has inertia degree $1$ in the extension $\QQ(\alpha_{\pp})$, where $\alpha_{\pp}$ is a root of the Frobenius polynomial $x^2-a_{\pp}x+\n(\pp)$, as required by the algorithm. At these primes, we have $a_{\pp_{2.1}}(E) = -1$ and    $a_{\pp_{3.1}}(E) =-2$, in agreement with the Hecke eigenvalues.
    
    \item For the relevant modulus $\nn_F = \pp_{2.1}^{5}\pp_{7.1}\pp_{7.2}\pp_{17.1}$, the ray class group is
      \[
      \CL(\nn_F) \simeq C_2 \times C_2 \times C_2 \times C_{24} \times C_{96}.
      \] 
    
    \item There are 31 index two subgroups of $\CL(\nn_F)$. For each of these and the full group, we consider the corresponding quadratic (or trivial) extension $L$ of $K$. We compute the modulus $\nn_L$  and the corresponding ray class group $\CL(\nn_L)$. For each of these ray class groups, we select a set of generators $\set{\chi_j}_{j=1}^{n}$ of the group of cubic characters. 
    For each of these generators, we pick a set of primes $\mathfrak{q}_j$ such that
    \[\langle \log(\chi_1(\mathfrak{q}_j)),..., \log(\chi_n(\mathfrak{q}_j)) \rangle_{j=1}^{n'} = (\ZZ/3\ZZ)^n \]
    for some $n'\in \ZZ$. In our example, we only needed to use the primes above $3$ and $23$. The eigenvalues for these primes are
\[
        a_{\pp_{3.1}}(F) = a_{\pp_{3.2}}(F) = -2, \quad a_{\pp_{23.1}}(F) = 0, \quad \text{and} \quad a_{\pp_{23.2}}(F) = 8,
\]
    which are all even, as required. Thus, we proceeded to the next step. 
    
    \item  For this step, we need a basis of quadratic characters $\set{\chi_i}_{i=1}^{n}$ of $\CL(\nn_F)$, together with a set of prime ideal $\set{\pp_i}$ coprime to $\nn_F$, such that 
    \[\set{\log(\chi_1(\mathfrak{q}_j),\log(\chi_2(\mathfrak{q}_j),\dots, \log(\chi_n(\mathfrak{q}_j) }_{j=1}^{2^n-1} = (\ZZ/2\ZZ)^n\backslash\set{0}. \]
    Here $n=5$, and it suffices to use the primes above the following rational primes:
\[
 \set{5, 11, 13, 19, 31, 53, 79, 89,  107, 131, 149, 157, 167, 257, 281, 457, 593, 1721}.
 \]
    For each of these, we check that the Hecke eigenvalues are equal to the trace of Frobenius of $E$, and we move on to the final step. 
    
    \item In the last step of algorithm, we need to check if the local $L$-factors of $E$ and $F$ are equal at primes dividing $2{\nn}(E){\nn}(F)\overline{{\nn}(F)}\Delta(K)$.  These are 
    \[
    \set{\pp_{2.1}, \pp_{7.1}, \pp_{7.2}, \pp_{17.1}}.
    \]
    We have already checked that the Hecke eigenvalues agree with the trace of Frobenius at all these except $\pp_{7.2}$, the unique prime dividing the level.    On the modular form side, the prime $\pp_{7.2}$ has Atkin Lehner eigenvalue $+1$. On the elliptic curve side, $E$ has non-split multiplicative reduction at $\pp_{7.2}$, and thus $a_{\pp_{7.2}}(E)=+1$, in agreement. So both $L$-functions have the same Euler factor at $\pp_{7.2}$.
\end{enumerate}
This completes the proof of the modularity of the elliptic curve $E$.

\begin{remark}
    As pointed out in the introduction, the results of \cite{caraiani2025modularityellipticcurvesimaginary} do not immediately imply that all elliptic curves over $\QQ(\sqrt{-17})$ are modular.  However, they do give sufficient conditions in terms of the mod-$3$ and mod-$5$ Galois representations, which are easily checked in individual cases.  For this elliptic curve $E$, the mod-$\ell$ Galois representations are surjective for all odd primes $\ell$ (see the curve's home page \lmfdbecnf{2.0.68.1}{7.2}{a}{2} for this information), which suffices, and gives a simpler way to establish modularity which was not available when this example was first studied.
\end{remark}
 \end{example}
 
\section{Results for other fields and future work}
\label{sec:otherfields}
In the previous section we have only given detailed results for one field, $\QQ(\sqrt{-17})$, whose class group is different from that of other imaginary quadratic fields for which any detailed computations have been published, to illustrate the general method.   However, our code works (at least in principle) over an arbitrary imaginary quadratic field, and we plan to produce more systematic data about newforms and their Hecke fields in the near future, data which we will publish via the~\lmfdb.  To date, we have at least some data for all $763$ fields of absolute discriminant up to $2500$, which includes all fields of class number up to $5$ except $\QQ(\sqrt{-2683})$,  with the largest class number being $63$ (for $\QQ(\sqrt{-2351})$) and the largest $2$-rank being $3$ (for 21 fields).   For each of these fields,  we have complete dimension data for all levels $\nn$ with norm up to a bound depending on the field, but we have only identified and computed eigenvalues for `rational' newforms $F$, that is, those with Hecke Field $k_F=\QQ$.  All this data for rational Bianchi newforms may be found in the \lmfdbbmftop.

The level norm bounds for the imaginary quadratic fields $K$ for which we have data are given in Table \ref{tab:level-ranges}, where $D=|\disc(K)|$.

\begin{table}[ht]
\caption{Level ranges for fields with $D<2500$.}
\label{tab:level-ranges}
\begin{tabular}{|C|CCCCCC|}
\hline
 D & 3 & 4 & 7,8,11 & 19,43 & 67 & 163   \\ \hline
 \n(\nn)\le{} & 150 000 & 100 000 & 50 000 & 15 000 & 10 000 & 5 000   \\
\hline 
\end{tabular}
\vskip10pt
\begin{tabular}{|C|CCCCCC|}
\hline
 D &  23 & 31 & 15-120 & 121-1000 & 1001-2100 & 2102-2500  \\ \hline
 \n(\nn)\le{} & 3 000 & 5 000 & 1 000 & 100 & 15 & 11  \\
\hline 
\end{tabular}
\end{table}


\begin{thebibliography}{{LMF}26}

\bibitem[AL78]{AtkinLi}
A.~O.~L. Atkin and W.~Li, \emph{Twists of newforms and pseudo-eigenvalues of
  {W}-operators}, Inventiones mathematicae \textbf{48} (1978), no.~3, 221--243.

\bibitem[Ara10]{aranes}
M.~Aranes, \emph{Modular symbols over number fields}, Ph.D. thesis, University
  of Warwick, December 2010, Available at
  \url{http://webcat.warwick.ac.uk/record=b2490253~S15}.

\bibitem[Ash77]{ash-def}
Avner Ash, \emph{Deformation retracts with lowest possible dimension of
  arithmetic quotients of self-adjoint homogeneous cones}, Math. Ann.
  \textbf{225} (1977), no.~1, 69--76. \MR{427490}

\bibitem[BBB{\etalchar{+}}21]{best2021computing}
Alex~J Best, Jonathan Bober, Andrew~R Booker, Edgar Costa, John~E Cremona,
  Maarten Derickx, Min Lee, David Lowry-Duda, David Roe, Andrew~V Sutherland,
  et~al., \emph{Computing classical modular forms}, Arithmetic {G}eometry,
  {N}umber {T}heory, and {C}omputation, Springer, 2021, pp.~131--213.

\bibitem[BCP97]{MAGMA}
Wieb Bosma, John Cannon, and Catherine Playoust, \emph{The {M}agma algebra
  system. {I}. {T}he user language}, J. Symbolic Comput. \textbf{24} (1997),
  no.~3-4, 235--265, Computational algebra and number theory (London, 1993).
  \MR{1 484 478}

\bibitem[BH07]{Berger}
Tobias Berger and Gergely Harcos, \emph{{$l$}-adic representations associated
  to modular forms over imaginary quadratic fields}, Int. Math. Res. Not.
  \textbf{23} (2007), rnm113.

\bibitem[Byg98]{bygott}
J.~S. Bygott, \emph{Modular forms and modular symbols}, Ph.D. thesis, The
  University of Exeter, 1998, Available at
  \url{http://hdl.handle.net/10871/8322}, p.~233.

\bibitem[CA14]{cremona_aranes}
J.~E. Cremona and M.~T. Aran\'{e}s, \emph{Congruence subgroups, cusps and
  {M}anin symbols over number fields}, Computations with modular forms,
  Contrib. Math. Comput. Sci., vol.~6, Springer, Cham, 2014, pp.~109--127.
  \MR{3381450}

\bibitem[CN25]{caraiani2025modularityellipticcurvesimaginary}
Ana Caraiani and James Newton, \emph{On the modularity of elliptic curves over
  imaginary quadratic fields}, \url{https://arxiv.org/abs/2301.10509}, 2025.

\bibitem[CPS20]{CremonaPageSutherland}
J.~E. Cremona, A.~Page, and A.~V. Sutherland, \emph{Sorting and labelling
  integral ideals in a number field},
  \url{https://doi.org/10.48550/arXiv.2005.09491}, 2020, Preprint.

\bibitem[Cre81]{cremonathesis}
J.~E. Cremona, \emph{Modular symbols}, Ph.D. thesis, The University of Oxford,
  1981, Available at
  \url{https://johncremona.github.io/theses/JCthesis-scan.pdf}.

\bibitem[Cre84]{cremona_hyptess}
\bysame, \emph{Hyperbolic tessellations, modular symbols, and elliptic curves
  over complex quadratic fields}, Compositio Math. \textbf{51} (1984), no.~3,
  275--324. \MR{743014}

\bibitem[Cre25]{Cremona_bianchi_progs}
\bysame, \emph{{bianchi-progs}},
  \url{https://github.com/JohnCremona/bianchi-progs}, 2025, GitHub repository
  of {\tt C++} code.

\bibitem[Cre26]{cremona-preprint}
\bysame, \emph{Hecke operators, {H}ecke eigensystems, and formal modular forms
  over number fields}, \url{https://arxiv.org/abs/2601.17524}, 2026.

\bibitem[\c{S}14]{sengun2014arithmetic}
Mehmet~Haluk \c{S}eng\"un, \emph{Arithmetic aspects of {B}ianchi groups},
  Computations with modular forms, Contrib. Math. Comput. Sci., vol.~6,
  Springer, Cham, 2014, pp.~279--315. \MR{3381457}

\bibitem[DGP10]{modularityexamplesDGP}
Luis Dieulefait, Lucio Guerberoff, and Ariel Pacetti, \emph{Proving modularity
  for a given elliptic curve over an imaginary quadratic field}, Math. Comp.
  \textbf{79} (2010), no.~270, 1145--1170. \MR{2600560}

\bibitem[EGM82]{elstrodtgrunewaldmennicke}
J.~Elstrodt, F.~Grunewald, and J.~Mennicke, \emph{On the group {${\rm
  PSL}_{2}({\bf Z}[i])$}}, Number {T}heory {D}ays, 1980 ({E}xeter, 1980),
  London Math. Soc. Lecture Note Ser., vol.~56, Cambridge Univ. Press,
  Cambridge-New York, 1982, pp.~255--283. \MR{697270}

\bibitem[GHM78]{grunewaldhellingmennicke}
F.~Grunewald, H.~Helling, and J.~Mennicke, \emph{{${\rm SL}_{2}$} over complex
  quadratic number fields. {I}}, Algebra i Logika \textbf{17} (1978), no.~5,
  512--580, 622. \MR{555260}

\bibitem[GL77]{GerardinLabesse}
P.~G\'erardin and J.-P. Labesse, \emph{The solution of a base change problem
  for {GL(2)} (following {L}anglands, {S}aito, {S}hintani)}, Automorphic Forms,
  Representations and L-functions (A.~Borel and W.~Casselman, eds.), Proceeding
  of {S}ymposia in {P}ure {M}athematics, vol. 33 (part 2), American
  Mathematical Society, 1977, pp.~115--133.

\bibitem[Gun99]{gunnells-qrank1}
Paul~E. Gunnells, \emph{Modular symbols for {${\bf Q}$}-rank one groups and
  {V}orono\u{\i} reduction}, J. Number Theory \textbf{75} (1999), no.~2,
  198--219. \MR{1681629}

\bibitem[HST93]{TaylorI}
Michael Harris, David Soudry, and Richard Taylor, \emph{{$l$}-adic
  representations associated to modular forms over imaginary quadratic fields.
  {I}. {L}ifting to {${\rm GSp}_4({\bf Q})$}}, Invent. Math. \textbf{112(2)}
  (1993), 377--411.

\bibitem[Kob84]{KoblitzECMF}
Neal Koblitz, \emph{Introduction to elliptic curves and modular forms},
  Graduate Texts in Mathematics, vol.~97, Springer-Verlag, New York, 1984.

\bibitem[Koe60]{koecher-cones}
Max Koecher, \emph{Beitr\"{a}ge zu einer {R}eduktionstheorie in
  {P}ositivit\"{a}tsbereichen. {I}}, Math. Ann. \textbf{141} (1960), 384--432.
  \MR{124527}

\bibitem[Kur78]{kurcanov}
P.~F. Kur\v{c}anov, \emph{The cohomology of discrete groups and {D}irichlet
  series that are related to {J}acquet-{L}anglands cusp forms}, Izv. Akad. Nauk
  SSSR Ser. Mat. \textbf{42} (1978), no.~3, 588--601. \MR{503433}

\bibitem[Lan76]{LangModularForms}
S.~Lang, \emph{Introduction to modular forms}, Springer-Verlag,
  Berlin--Heidelberg--New York, 1976, Grundlehren der mathematischen
  Wissenschaften 222.

\bibitem[Lin05]{lingham}
Mark~Peter Lingham, \emph{Modular forms and elliptic curves over imaginary
  quadratic fields}, Ph.D. thesis, The University of Nottingham, 2005,
  Available at \url{http://eprints.nottingham.ac.uk/10138/}, p.~120.

\bibitem[{LMF}26]{lmfdb}
The {LMFDB Collaboration}, \emph{The {L}-functions and {M}odular {F}orms
  {D}atabase}, \url{https://www.lmfdb.org}, 2026, [Online; accessed 7 January
  2026].

\bibitem[Miy71]{Miyake}
Toshitsune Miyake, \emph{On automorphic forms on $\mathrm{GL}_2$ and {H}ecke
  operators}, Annals of Mathematics \textbf{\textbf{94}} (1971), no.~1,
  174--189.

\bibitem[Ser73]{serre-book}
J.-P. Serre, \emph{A course in arithmetic}, Graduate Texts in Mathematics, No.
  7, Springer-Verlag, New York-Heidelberg, 1973, Translated from the French.
  \MR{0344216}

\bibitem[Swa71]{swan}
Richard~G. Swan, \emph{Generators and relations for certain special linear
  groups}, Advances in Math. \textbf{6} (1971), 1--77. \MR{284516}

\bibitem[Tay94]{TaylorII}
Richard Taylor, \emph{{$l$}-adic representations associated to modular forms
  over imaginary quadratic fields. {II}}, Invent. Math. \textbf{116(1-3)}
  (1994), 619--643.

\bibitem[Tha23]{kalanithesis}
Kalani Thalagoda, \emph{Computational aspects of {B}ianchi modular forms},
  Ph.D. thesis, The University of North Carolina at Greensboro, 2023, Avalable
  at \url{https://libres.uncg.edu/ir/uncg/listing.aspx?id=47004}.

\bibitem[Wei71]{Weil}
A.~Weil, \emph{Dirichlet series and automorphic forms}, Lecture Notes in
  Mathematics, no. 189, Springer-Verlag, 1971.

\bibitem[Whi90]{whitley}
Elise Whitley, \emph{Modular forms and elliptic curves over imaginary quadratic
  number feilds}, Ph.D. thesis, The University of Exeter, 1990, Available at
  \url{http://hdl.handle.net/10871/8427}, p.~233.

\bibitem[Yas10]{yasak_hyptess}
Dan Yasaki, \emph{Hyperbolic tessellations associated to {B}ianchi groups},
  Algorithmic number theory, Lecture Notes in Comput. Sci., vol. 6197,
  Springer, Berlin, 2010, pp.~385--396. \MR{2721434}

\end{thebibliography}

\newcommand{\etalchar}[1]{$^{#1}$}
\providecommand{\bysame}{\leavevmode\hbox to3em{\hrulefill}\thinspace}
\providecommand{\MR}{\relax\ifhmode\unskip\space\fi MR }
\providecommand{\MRhref}[2]{%
  \href{http://www.ams.org/mathscinet-getitem?mr=#1}{#2}
}
\providecommand{\href}[2]{#2}

\end{document}